\documentclass[aap,MSNbibl,seceqn,dvips]{arximspdf}
\usepackage{graphicx}
\doi{10.1214/13-AAP952} 
\volume{24}
\issue{4}
\pubyear{2014}
\firstpage{1446}
\lastpage{1481}

\makeatletter

\newcommand{\rrvert}{\vert}
\newcommand{\llvert}{\vert}
\newtheorem{theorem}{Theorem}[section]
\newtheorem{lemma}[theorem]{Lemma}
\newtheorem{proposition}[theorem]{Proposition}
\newtheorem{corollary}[theorem]{Corollary}
\newcommand{\eqref}[1]{(\ref{#1})}
\newtheorem{lemmas}{Lemma}[section]
\def\p{\mathbb{P}}

\def\ee{\mathrm{e}}
\def\d{\mathrm{d}}

\def\R{\mathbb{R}}
\def\N{\mathbb{N}}
\def\E{\mathbb{E}}
\def\equilaw{{\stackrel{\mathrm{law}}{=}}}
\makeatother

\begin{document}
\begin{frontmatter}

\title{Poisson--Dirichlet statistics for the extremes of a~log-correlated Gaussian field}
\runtitle{Poisson--Dirichlet statistics, log-correlated Gaussian field}

\begin{aug}
\author[A]{\fnms{Louis-Pierre} \snm{Arguin}\corref{}\thanksref{t1}\ead[label=e1]{arguinlp@dms.umontreal.ca}}
\and
\author[B]{\fnms{Olivier} \snm{Zindy}\thanksref{t2}\ead[label=e2]{olivier.zindy@upmc.fr}}
\thankstext{t1}{Supported by a NSERC discovery grant and a grant FQRNT
\textit{Nouveaux chercheurs}.}
\thankstext{t2}{Supported in part by the French ANR project MEMEMO2
2010 BLAN 0125.}
\runauthor{L.-P. Arguin and O. Zindy}
\affiliation{Universit\'e de Montr\'eal and Universit\'e Paris 6}
\address[A]{D\'epartement de Math\'ematiques\\
\quad et Statistique\\
Universit\'e de Montr\'eal\\
Montr\'eal, Quebec H3T 1J4\\
Canada\\
\printead{e1}} 
\address[B]{Laboratoire de Probabilit\'es\\
\quad et Mod\`eles Al\'eatoires\\
CNRS UMR 7599\\
Universit\'e Paris 6\\
4 place Jussieu\\
75252 Paris Cedex 05\\
France\\
\printead{e2}}
\end{aug}

\received{\smonth{6} \syear{2013}}

%
\begin{abstract}
We study the statistics of the extremes of a discrete Gaussian field
with logarithmic correlations at the level of the Gibbs measure.
The model is defined on the periodic interval $[0,1]$, and its
correlation structure is nonhierarchical.
It is based on a model introduced by
Bacry and Muzy [\textit{Comm. Math. Phys.} \textbf{236} (2003) 449--475]
(see also Barral and Mandelbrot [\textit{Probab. Theory Related Fields}
\textbf{124} (2002) 409--430]), and is similar to the logarithmic Random Energy Model
studied by Carpentier and Le Doussal
[\textit{Phys. Rev. E} (3) \textbf{63} (2001) 026110]
and more recently by Fyodorov and Bouchaud
[\textit{J. Phys. A} \textbf{41} (2008) 372001].
At low temperature, it is shown that the normalized covariance of two
points sampled from the Gibbs measure is either $0$ or $1$.
This is used to prove that the joint distribution of the Gibbs weights
converges in a suitable sense to that of a Poisson--Dirichlet variable.
In particular, this proves a conjecture of Carpentier and Le Doussal
that the statistics of the extremes of the log-correlated field behave
as those of i.i.d. Gaussian variables and of branching Brownian motion
at the level of the Gibbs measure.
The method of proof is robust and is adaptable to other log-correlated
Gaussian fields.
\end{abstract}

%

\end{frontmatter}

\section{Introduction}\label{sec1}

This paper studies the statistics of the extremes of a Gaussian field
whose correlations decay logarithmically with the distance.
The model is related to the process introduced by Bacry and Muzy \cite
{bacry-muzy} (see also Barral and Mandelbrot \cite{barral-mandelbrot})
and is similar to the \textit{logarithmic random energy model} or \textit
{log-REM} studied by Carpentier and Le Doussal \cite
{carpentier-ledoussal}, and Fyodorov and Bouchaud \cite{bouchaud-fyodorov}.
Another important log-correlated model is the two-dimensional discrete
Gaussian free field.

The statistics of the extremes of log-correlated Gaussian fields are
expected to resemble those of i.i.d. Gaussian variables or \textit
{random energy model} (REM) and at a finer level, those of branching
Brownian motion.
In fact, log-correlated fields are conjectured to be the critical case
where correlations start to affect the statistics of the extremes.
The reader is referred to the works of Carpentier and Le Doussal \cite
{carpentier-ledoussal};
Fyodorov and Bouchaud \cite{bouchaud-fyodorov}; and Fyodorov, Le
Doussal and Rosso~\cite{fyodorov-ledoussal-rosso} for physical
motivations of this fact.
The analysis for general log-correlated Gaussian field is complicated
by the fact that, unlike branching Brownian motion, the correlations do
not necessarily exhibit a tree structure.

The approach of this paper is in the spirit of the seminal work of
Derrida and Spohn \cite{derrida-spohn} who studied the extremes of
branching Brownian motion using the Gibbs measure.
The method of proof presented here is robust and applicable to a large
class of nonhierarchical log-correlated fields.
The model studied here has the advantages of having a graphical
representation of the correlations, a continuous scale parameter and no
boundary effects (cf. Section~\ref{sec1.1}) which make the ideas of the method
more transparent.
Even though correlations are not tree-like for general log-correlated
models, such fields can often be decomposed as a sum of independent
fields acting on different scales.
The main results of the paper are Theorem \ref{thm:low-temperature} on
the correlations of the extremes and Theorem \ref{thm:main4} on the
statistics of the Gibbs weights.
The results show that, in effect, the statistics of the extremes of the
log-correlated field are the same as those of branching Brownian motion
at the level of the Gibbs measure,
as conjectured by Carpentier and Le Doussal~\cite{carpentier-ledoussal}.

The method of proof is outlined in Section~\ref{sect:outline}.
The proof of the first theorem is based on an adaptation of a technique
of Bovier and Kurkova \cite{bovier-kurkova1,bovier-kurkova2}
originally developed for hierarchical Gaussian fields such as branching
Brownian motion.
For this purpose, we need to introduce a family of log-correlated
Gaussian models where the variance of the fields in the
scale-decomposition depends on the scale.
The free energy of the perturbed models is computed using ideas of
Daviaud \cite{daviaud}.
The second theorem on the Poisson--Dirichlet statistics of the Gibbs
weights is proved using the first theorem on correlations and general
spin glass theory results.

\subsection{A log-correlated Gaussian field}\label{sec1.1}

Following \cite{bacry-muzy}, we consider the half-infinite cylinder
\[
\mathcal C^+:=\bigl\{(x,y); x \in[0,1]_\sim, y \in\R_+^* \bigr\},
\]
where $[0,1]_\sim$ stands for the unit interval where the two endpoints
are identified.
We write $\|x-x'\|:=\min\{|x-x'|, 1-|x-x'|\}$ for the distance on
$[0,1]_\sim$.

The following measure is put on $\mathcal C^+$:
\[
\theta(\mathrm{d} x,\mathrm{d} y):= y^{-2} \,\d x \,\d y.
\]
For $\sigma>0$, \textit{the variance parameter}, there exists a random
measure $\mu$ on $\mathcal C^+$ that satisfies:
\begin{longlist}[(ii)]
\item[(i)] for any measurable set $A$ in $\mathcal B (\mathcal C^+)$,
the random variable $\mu(A)$ is a centered Gaussian with variance
$\sigma^2 \theta(A)$;
\item[(ii)] for every sequence of disjoint sets $(A_n)_n$ in $\mathcal
B (\mathcal C^+)$, the Borel $\sigma$-algebra associated with
$\mathcal
C^+$, the random variables $(\mu(A_n))_n$ are independent and
\[
\mu\biggl(\bigcup_{n}
A_n \biggr)= \sum_{n}
\mu(A_n)\qquad \mbox{a.s.}
\]
\end{longlist}
Let $\Omega$ be the probability space on which $\mu$ is defined, and
let $\p$ be the law of $\mu$.
The space $\Omega$ is endowed with the $\sigma$-algebras $\mathcal F_u$
generated by the random variables $\mu(A)$, for all the sets $A$ at a
distance greater than $u$ from the $x$-axis.
The reader is referred to \cite{bacry-muzy} for the existence of the
probability space $(\Omega, (\mathcal F_u)_u, \p)$.

\begin{figure}

\includegraphics{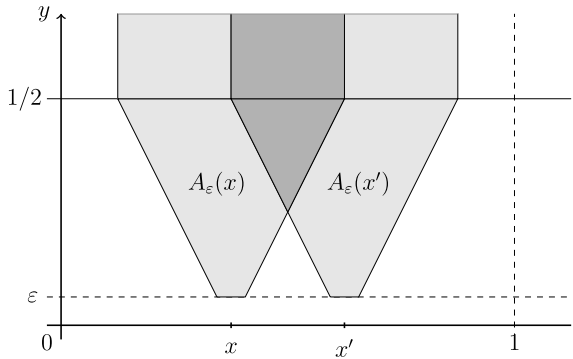}

\caption{The two subsets $A_\varepsilon(x)$ and $A_\varepsilon(x')$ for
$\varepsilon=1/N$.
The variance of the variables is given by the integral over $\theta(\d
t,\d y)=y^{-2}\,\d t \,\d y$ of the lighter gray area above $\varepsilon
=1/N$, and
the covariance by the integral over the intersection of the subsets,
the darker gray region.}
\label{fig: cones}
\end{figure}

The subsets needed for the definition of the Gaussian field are the
cone-like subsets $A_u (x)$ of $\mathcal C^+$,
\[
A_u (x):=\bigl\{(s,y) \in\mathcal C^+ \dvtx y \ge u, -f(y)/2 \le
s-x \le f(y)/2 \bigr\},
\]
where $f(y)= y$ for $y \in(0,1/2)$ and $f(y)= 1/2$ otherwise. See
Figure~\ref{fig: cones} for a depiction of the subsets.
Observe that, by construction, if $\|x-x'\|=\ell>u$, then $A_u (x)$ and
$A_u(x')$ intersect exactly above the line $y=\ell$.

The Gaussian process $\omega_u= (\omega_u(x), x \in[0,1]_\sim)$
is defined using the random measure $\mu$,
%
%
\begin{equation}
\label{eqn: omega} \omega_u(x):= \mu\bigl(A_u (x)\bigr),\qquad  x
\in[0,1]_\sim.
\end{equation}
By properties (i) and (ii) of $\mu$ listed above, the covariance
between $\omega_u(x)$ and $\omega_u(x')$ is given by the integral over
$\theta$ of the intersection of $A_u(x)$ and $A_u(x')$,
%
%
\begin{equation}
\label{eqn: corr} \E\bigl[ \omega_u(x)\omega_u
\bigl(x'\bigr)\bigr] = \int_{A_u(x)\cap A_u(x')} \theta(\d s,\,\d
y).
\end{equation}


%


The paper focuses on a discrete version of $\omega_u$.
Let $N\in\N$, and take $\varepsilon=1/N$. Define the set
\[
\mathcal X_{N}=\mathcal X_{\varepsilon}:= \biggl\{0,
\frac{1}{N},\frac
{2}{N},\ldots,\frac{i}{N},\ldots,
\frac{N-1}{N} \biggr\}.
\]
The notation $\mathcal X_{N}$ and $\mathcal X_{\varepsilon}$ will be
used equally depending on the context.
For a given~$N$, the $\log$-correlated Gaussian field is the collection
of Gaussian centered random variables $\omega_{\varepsilon}(x)$ for
$x\in\mathcal{X}_N$,
%
%
\begin{equation}
\label{eqn: X} X=(X_x, x\in\mathcal X_N)=\bigl(
\omega_\varepsilon(x), x \in\mathcal X_{N}\bigr).
\end{equation}
A compelling feature of this construction is that a \textit{scale
decomposition} for $X$ is easily obtained from property (ii) above.
Indeed, it suffices to write the variable~$X_x$ as a sum of independent
Gaussian fields corresponding to disjoint horizontal strips of
$\mathcal C^+$.
The $y$-axis then plays the role of the scale.

The covariances of the field are computed from \eqref{eqn: corr} by
straightforward integration; see also Figure~\ref{fig: cones}.
%
%
\begin{lemma}
\label{lem:covariance}
For any $0<\varepsilon=1/N<1/2$,
\begin{eqnarray*}
\E\bigl[X_x^2\bigr] &=& \sigma^2 ( \log N
+1-\log2),\qquad x \in\mathcal X_N,
\\
\E[ X_xX_{x'} ] &=& \sigma^2\bigl( \log
\bigl(1/ \bigl\Vert x-x' \bigr\Vert\bigr)-\log2\bigr),\qquad  x \neq x'
\in\mathcal X_N.
\end{eqnarray*}
\end{lemma}
Similar constructions of log-correlated Gaussian fields using a random
measure on cone-like subsets are also possible in two dimensions; see,
for example, \cite{robert-vargas}.

\subsection{Main results}\label{sec1.2}
Without loss of generality, the results of this section are stated for
the variance parameter $\sigma=1$.
The points where the field is unusually high, \textit{the extremes} or
\textit{the high points}, can be studied using a minor adaptation of
the arguments of Daviaud for the
two-dimensional discrete Gaussian free field~\cite{daviaud}. We denote
by $ \vert\mathcal{A}\vert$ the cardinality of a finite set $
\mathcal{A}$.
%
%
\begin{theorem}[(Daviaud \cite{daviaud})]
\label{thm:main2}
Let
\[
\mathcal{H}_N(\gamma):= \{ x \in\mathcal{X}_N \dvtx
X_x \ge\sqrt{2} \gamma\log N \}
\]
be the set of $\gamma$-high points. Then for any $0<\gamma<1$,
\[
\lim_{N \to\infty} \frac{\log\vert\mathcal{H}_N(\gamma)\vert
}{\log
N}=1-\gamma^2 \qquad\mbox{in probability.}
\]
Moreover, for all $\rho>0$ there exists a constant $c=c(\rho)>0$ such that
\[
\p\bigl(\bigl\vert\mathcal{H}_N(\gamma)\bigr\vert\le N^{(1-\gamma^2)-\rho}
\bigr) \le\exp\bigl\{- c (\log N)^2\bigr\}
\]
for $N$ large enough.
\end{theorem}
The technique of Daviaud is based on a tree approximation introduced by
Bolthausen, Deuschel and Giacomin \cite{bolthausen-deuschel-giacomin}
for the discrete two-dimensional Gaussian free field.
There, the technique is used to obtain the first order of the maximum.
The same argument applies here. Theorem \ref{thm:main2} and simple
Gaussian estimates yield
%
%
\begin{equation}
\label{eqn: first order} \lim_{N\to\infty} \frac{\max_{x \in
\mathcal X_N} X_x}{\log N}= \sqrt{2}\qquad \mbox{a.s.}
\end{equation}
The important feature of Theorem \ref{thm:main2} and equation \eqref
{eqn: first order} is that they
are identical to the results for $N$ i.i.d. Gaussian variables of
variance $\log N$.
In other words, the above observables of the high points are not
affected by the correlations of the field.
The i.i.d. case is called the \textit{random energy model} (REM) in
the spin glass literature.

The starting point of the paper is to understand to which extent
i.i.d. statistics is a good approximation for more refined observables
of the extremes of log-correlated Gaussian fields.
To this end, we turn to tools of statistical physics which allow for a
good control of the correlations.

First, consider the \textit{partition function} $Z_N(\beta)$ of the
model ($\beta$ stands for the inverse-temperature),
\[
Z_N(\beta):= \sum_{x \in\mathcal X_N} \exp\{ \beta
X_x \}\qquad \forall\beta>0,
\]
and the \textit{free energy}
\[
f_N(\beta):= \frac{1}{\log N} \log Z_N(\beta)\qquad
\forall\beta>0.
\]
Theorem \ref{thm:main2} is used to compute the free energy of the model.
%
%
\begin{corollary}
\label{cor:main3}
Let $\beta_c:=\sqrt{2}$. Then, for all $\beta>0$
\[
f(\beta):= \lim_{N \to\infty} f_N(\beta)=
\cases{ \displaystyle 1 +
\frac{\beta^2}{2}, &\quad $\mbox{if } \beta< \beta_c$, \vspace*{2pt}
\cr
\sqrt{2} \beta, & \quad $\mbox{if } \beta\ge\beta_c$,}\qquad \mbox{a.s. and in
$L^1$.}
\]
\end{corollary}
The free energy is the same as for the REM with variance $\log N$.
In particular, the model undergoes \textit{freezing} above $\beta_c$ in
the sense that the quantity $f(\beta)/\beta$ is constant.

More importantly, consider the \textit{normalized Gibbs weights} or
\textit{Gibbs measure}
\[
G_{\beta, N}(x):=\frac{\ee^{\beta X_x}}{Z_N(\beta)},\qquad x \in\mathcal X_N.
\]
By design, the Gibbs measure concentrates on the high points of the
Gaussian field.
The first main result of the paper is to achieve a control of the
correlations at the level of the Gibbs measure.
Precisely, with spin glasses in mind, we consider the normalized
covariance or \textit{overlap}
%
%
\begin{equation}
\label{eqn: q} q(x,y)=q^{(N)}(x,y): = - \frac{\log\Vert y -x \Vert}
{\log N}, \qquad x,y \in
\mathcal X_N.
\end{equation}
Clearly, $\|x-y\|=\varepsilon^{q(x,y)}$ and $0\leq q(x,y) \leq1$.
Moreover, the overlap $q(x,y)$ is equal to the normalized correlations
$\E[ X_x X_y] / \E[ X_x^2] $ plus a term that goes to zero as $N$ goes
to infinity.

A fundamental object, that records the correlations of high points, is
the \textit{distribution function of the overlap} sampled from the
Gibbs measure.
Namely, denote by $G_{\beta, N}^{\times2}$ the product measure on
$\mathcal X_N\times\mathcal X_N$. Let $(x_1,x_2)\in\mathcal X_N^2$ be
sampled from $G_{\beta, N}^{\times2}$. Write for simplicity $q_{12}$
for $q(x_1,x_2)$. The averaged distribution function of the overlap is
%
%
\begin{equation}
\label{eqn: x} x_{\beta}^{(N)}(q):= {\mathbb E}\bigl[
G_{\beta, N}^{\times2} \{q_{12} \le q \} \bigr], \qquad 0\leq q
\leq1.
\end{equation}

The first result is the analogue of results of Derrida and Spohn for
the Gibbs measure of branching Brownian motion (see equation (6.19) in
\cite{derrida-spohn}),
of Chauvin and Rouault on branching random walks \cite{chauvin-rouault}
and of Bovier and Kurkova on Derrida's
\textit{generalized random energy models} (GREM) \cite{derrida,bovier-kurkova1}.
It had been conjectured for nonhierarchical log-correlated Gaussian
field by Carpentier and Le Doussal; see page 16 in \cite{carpentier-ledoussal}.
%
%
\begin{theorem}
\label{thm:low-temperature}
For $\beta> \beta_c$,
\[
\lim_{N\to\infty}x_{\beta}^{(N)}(q)=\lim
_{N\to\infty} {\mathbb E}\bigl[ G_{\beta, N}^{\times2}
\{q_{12} \le q \} \bigr]= \cases{ \displaystyle\frac{\beta_c}{\beta},&\quad
$\mbox{for $0
\le q <1$,}$\vspace*{2pt}
\cr
1, &\quad $\mbox{for $q=1$.}$}
\]
\end{theorem}
This result is the same as for the REM model \cite{talagrand}.
It is therefore consistent with rich statistics of extremes consisting
of many high values order one away of each other and whose correlations
are either very high or close to $0$.
This result is in expectation.
The typical behavior of the random variable $G_{\beta,N}^{{ \times
2}}\{
q_{12}\leq q\}$ for $q$ small in terms of $\beta$ should be
exponentially small in $\beta$ rather than $1/\beta$.
To see this, at the heuristic level, it is informative to consider the
i.i.d. case where the same phenomenon occurs.
Consider $N$ i.i.d. Gaussian random variables $(X_i)_{1 \le i \le N}$
of variance $\log N$ ordered in a decreasing way.
In this case, $q_{ij}=0$ if $i\neq j$. The following inequality is
easily verified:
\[
G_{\beta,N}^{\times2}\{q_{12}=0\}=\frac{\sum_{i\neq j} \mathrm{e}^{\beta
X_i}\mathrm{e}^{\beta X_j}}{ (\sum_i \mathrm{e}^{\beta X_i} )^2}\leq2
\sum_{j\geq2}\mathrm{e}^{\beta(X_j-X_1)}.
\]
In particular, since the gap $X_1-X_2$ is of order one in the limit and
since the density of points at distance $x$ from the maximum is bounded
by $\mathrm{e}^{Cx}$ for $C$ large enough (see \cite{bolthausen-sznitman} for a
precise statement in terms of extremal process), the typical behavior
of $G_{\beta,N}^{{ \times2}}\{q_{12}=0\}$ is expected to be
exponentially small in $\beta$.

We remark also that for $\beta\leq\beta_c$ the free energy contains
all information about the two-overlap distribution. Indeed, since the
free energy in Corollary \ref{cor:main3} is differentiable for every
$\beta>0$ including $\beta_c$, we have by the convexity of the free
energy that the derivative of the limit is the limit of the
derivatives. Hence
\[
\lim_{N\to\infty} f'_N(\beta)=\lim
_{N\to\infty}\beta\bigl(1-\E G_{\beta
,N}^{\times2}[q_{12}]
\bigr)=f'(\beta).
\]
The first equality is by Gaussian integration by part.
It follows that\break  $ \lim_{N} \E[G_{\beta,N}^{\times2}(q_{12})
]=0$ for $\beta\leq\beta_c$.
In particular, since the correlations are positive, the overlap of two
sampled points is $0$ almost surely for every $\beta\leq\beta_c$.

In the case of $\beta>\beta_c$, the first moment of the two-overlap
distribution is strictly greater than $0$, therefore more information
is needed to determine the distribution.
One way to proceed would be to obtain enough expectations of functions
of $q_{12}$ to determine the distribution.
This can be done by adding parameters to the field and consider the
appropriate derivative of the free energy of the perturbed model.
This is similar in spirit to the $p$-spin perturbations for the
Sherrington--Kirkpatrick model in spin glasses; see, for example, \cite
{talagrand}.
It turns out that this kind of pertubative approach pioneered by Bovier
and Kurkova in \cite{bovier-kurkova2} for Gaussian fields on trees can
be generalized to log-correlated fields.
The control of the correlations is achieved by introducing a perturbed
version of the model at a specific scale; cf. Section~\ref{sect: perturbed}.
In the present case, the proof is more intricate since the structure of
correlations of the Gaussian field for finite $N$ is not tree-like or
\textit{ultrametric} as in the cases of branching Brownian motion and GREM's.
For example, for branching Brownian motion, $q(x,y)$ corresponds to the
branching time of the common ancestor of two particles at time $t$, $x$
and $y$, divided by $t$. Because of the branching structure,
%
%
\begin{equation}
\label{eqn: ultra} \qquad\mbox{the inequality } q(x,y) \geq\min\bigl\{
q(x,z),q(y,z)\bigr\}
\mbox{ is satisfied for all $x,y,z$.}
\end{equation}
[The terminology \textit{ultrametric} comes from the fact that the
distance induced by the form $q(\cdot,\cdot)$ is ultrametric when
\eqref
{eqn: ultra} holds.]

The \textit{Parisi ultrametricity conjecture} in the spin-glass
literature states that, even though tree-like correlations might not be
present for finite $N$, ultrametric correlations are recovered in the
limit $N\to\infty$ for a large class of Gaussian fields at the level of
the Gibbs measure, that is,
%
%
\begin{equation}
\label{eqn: ultra+} \lim_{N\to\infty} {\mathbb E}\bigl[ G_{\beta,
N}^{\times3}
\bigl\{q_{12} \geq\min\{q_{13},q_{23}\} \bigr\}
\bigr] =1.
\end{equation}
It is not hard to see that Theorem \ref{thm:low-temperature} implies
the ultrametricity conjecture for the Gaussian field considered, since
the overlaps can only take value $0$ or $1$.
(In the language of spin glasses, the field is said to admit a \textit
{one-step replica symmetry breaking} at low temperature.)

The second main result describes the joint distribution of overlaps
sampled from the Gibbs measure.
To this end, for $s\geq2$, we denote the product of Gibbs measure on
$\mathcal X_N^s$ by $G_{\beta,N}^{\times s}$.
We consider the class of continuous functions $F\dvtx[0,1]^{
{s(s-1)}/{2}}\to\R$.
We write $\E G_{\beta,N}^{\times s} [F(q_{ll'}) ]$ for $\E
G_{\beta,N}^{\times s} [F(\{q(x_l,x_{l'})\}_{1 \le l<l' \le s})
]$, that is, the averaged expectation of $F(\{q(x_l,x_{l'})\}_{1 \le l<l'
\le s})$ when $(x_1,\ldots,x_s)$ is sampled from $G_{\beta,N}^{\times s}$.
We recall the definition of a Poisson--Dirichlet variable.
For $0<\alpha<1$, let $\eta=(\eta_i,i\in\N)$ be the atoms of a Poisson
random measure on $(0,\infty)$ of intensity measure $s^{-\alpha-1} \,
\mathrm{d}s$.
A \textit{Poisson--Dirichlet variable} $\xi$ of parameter $\alpha$ is a
random variable on the space of decreasing weights $\vec
{s}=(s_1,s_2,\ldots)$ with $1\geq s_1\geq s_2\geq\cdots\geq0$ and
$\sum_{i}s_i\leq1$ which has the same law as
\[
\xi\,\equilaw\,\biggl(\frac{\eta_i}{\sum_j\eta_j}, i\in\N\biggr
)_\downarrow,
\]
where $\downarrow$ stands for the decreasing rearrangement.
%
%
\begin{theorem}
\label{thm:main4}
Let $\beta>\beta_c$ and $\xi=(\xi_k,k\in\N)$ be a Poisson--Dirichlet
variable of parameter $\beta_c/\beta$.
Denote by $E$ the expectation with respect to $\xi$.
For any continuous function $F\dvtx[0,1]^{{s(s-1)}/{2}}\to\R$
of the
overlaps of $s$ points,
\[
\lim_{N\to\infty} {\mathbb E}G_{\beta, N}^{\times s} \bigl[
F(q_{ll'}) \bigr] = E \biggl[ \sum_{k_1\in\N,\ldots,k_s\in\N}
\xi_{k_1}\cdots\xi_{k_s} F(\delta_{k_lk_{l'}}) \biggr].
\]
\end{theorem}
It is important to stress that, as in the case of branching Brownian
motion and unlike the REM,
it is not the collection $(G_{\beta,N}(x), x \in\mathcal X_N
)_\downarrow$ \textit{per se} that converges to a Poisson--Dirichlet variable.
Rather, the result suggests that the Poisson--Dirichlet weights are
formed by the sum of the Gibbs weights of high points that are
arbitrarily close to each other because the continuity of the function
$F$ naturally identifies points $x,y$ for which $q(x,y)$ tends to $1$
in the limit $N\to\infty$.
In the theory of spin glasses, these clusters of high points are often
called \textit{pure states}.\vadjust{\goodbreak}
For more on the connection with spin glasses, the reader is referred to
\cite{talagrand-pure} where the pure states are constructed explicitly
for mean-field models.

\subsection{Relation to previous results}\label{sec1.3}

Bolthausen and Kistler have studied a family of models called \textit
{generalized GREMs} for which the correlations are not ultrametric
\cite
{bolthausen-kistler1,bolthausen-kistler2} for finite $N$.
By construction, the overlaps of these models can only take a finite
number of values (uniformly in $N$, the number of variables). They
compute the free energies and the Gibbs measure and prove the
Parisi ultrametricity conjecture for these.
Bovier and Kurkova \cite{bovier-kurkova1,bovier-kurkova2} have obtained
the distribution of the Gibbs measure for Gaussian fields, called the
CREMs, where the values of the overlaps are not a priori restricted.
Their analysis is restricted to models with ultrametric correlations
and include the case of branching Brownian motion.

The works of Bolthausen, Deuschel and Zeitouni \cite
{bolthausen-deuschel-zeitouni}, Bramson and Zei\-touni~\cite
{bramson-zeitouni} and Ding \cite{ding}
establish the tightness of the recentered maximum of the
two-dimensional discrete Gaussian free field. We expect that their
method can be applied to the Gaussian field we consider.

We note that Fang and Zeitouni \cite{fang-zeitouni} have studied a
branching random walk model where the variance of the motion is time-dependent.
This model is related to the simpler GREM model of spin glasses and to
the CREM of Bovier and Kurkova. The family of log-correlated Gaussian
fields introduced in Section~\ref{sec2.2} is akin to these hierarchical models,
where the scale parameter replaces the time parameter.

\section{Outline of the proof}
\label{sect:outline}

The proof is split in three steps, and each can be adapted (with
different correlation estimates) to other log-correlated Gaussian fields.
The Gaussian field we study has a graphical representation of its
correlations as well as no boundary effect which help in illustrating
the method.
\subsection{A family of perturbed models}
\label{sect: perturbed}
In this section, we define a family of Gaussian fields for which the
variance parameter $\sigma$ is scale-dependent.
It can be seen as the GREM analogue for the nonhierarchical Gaussian
field considered here.
We restrict ourselves to the case where $\sigma$ takes two values,
which is the one needed for the proof of Theorem \ref{thm:low-temperature}.
However, the construction and the results can hold for any finite
number of values.

Fix $\varepsilon= 1/N$. We introduce a scale (or time) parameter $t$ by
defining for any $t \in[0,1]$,
\[
X_x(t):= \omega_{\varepsilon^t}(x),\qquad x\in\mathcal{X}_\varepsilon.
\]
Observe that for any fixed $x$, the process $(X_x(t))_{0 \le t \le1}$
has independent increments and is a martingale for the filtration
$(\mathcal F_{\varepsilon^t}, t\geq0)$,
\[
\E\bigl[ X_x(t) | \mathcal F_{\varepsilon^s}\bigr]=
X_x(s)\qquad \mbox{for $t>s$}.
\]
This is a consequence of the defining property (ii) of the random
measure~$\mu$.\vadjust{\goodbreak}

\begin{figure}

\includegraphics{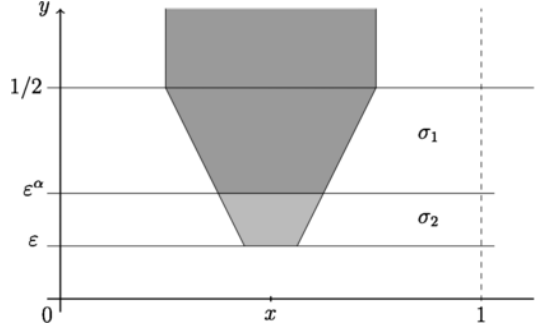}

\caption{The cone associated with the process $ Y_x(\cdot)$.}
\label{fig:perturbed-daviaud1}
\end{figure}

The parameters of the family of perturbed models are $\alpha$ where $0
< \alpha<1$
and $\vec{\sigma}=(\sigma_1,\sigma_2)$ with $\sigma_i >0$, $i=1,2$.
For the sake of clarity and to avoid repetitive trivial corrections, it
is assumed throughout the paper that $N^{\alpha}$ and $N^{1-\alpha}$
are integers.
The Gaussian field $Y^{(\vec{\sigma},\alpha)}(t)=(Y^{(\vec{\sigma
},\alpha)}_x(t), x\in\mathcal X_\varepsilon)$ is defined from the
field $X$ as follows:
%
%
\begin{equation}
\label{eqn: perturbed Y} Y^{(\vec{\sigma},\alpha)}_x(t)=\cases{ %
\sigma_1 X_x(t), &\quad $\mbox{if } 0 < t
\le\alpha,$
\vspace*{2pt}\cr
\sigma_1 X_x(\alpha) + \sigma_2 \bigl(
X_x(t)-X_x(\alpha) \bigr), &\quad $\mbox{if } \alpha< t \le
1.$}
\end{equation}
The construction is depicted in Figure~\ref{fig:perturbed-daviaud1}.
We write $Y^{(\vec{\sigma}, \alpha)}$ for the field $(Y^{(\vec
{\sigma},
\alpha)}_x(1),  x\in\mathcal X_\varepsilon)$.
The dependence on $\vec{\sigma}$ and $ \alpha$ will sometimes be
dropped in the notation of $Y$ for simplicity.


%


Consider the partition function $Z^{(\vec{\sigma}, \alpha)}_N(\beta)$
of the perturbed model
%
%
\begin{equation}
\label{eq:partition-function-perturbed} Z_N^{(\vec{\sigma}, \alpha
)}(\beta):=\sum
_{x\in\mathcal X_N} \exp(\beta Y_x ),
\end{equation}
and the free energy
\[
f^{(\vec{\sigma}, \alpha)}_N(\beta):= \frac{1}{\log N} \log
Z_N^{(\vec
{\sigma}, \alpha)}(\beta)\qquad \forall\beta>0.
\]

The log number of high points can be computed for the Gaussian field
$Y$ using Daviaud's technique recursively.
The free energy is then obtained by doing an explicit sum on these high points.
This is the object of Sections~\ref{sec:perturbed-maximum} and~\ref
{sect: free energy}.
The result is better expressed in terms of the free energy of the REM
with $N$ i.i.d. Gaussian variables of variance $\sigma^2\log N$,
\[
f\bigl(\beta; \sigma^2\bigr):= \cases{\displaystyle 1+\frac{\beta^2 \sigma
^2}{2}, &\quad$\mbox{if $\beta\leq\beta_c\bigl(\sigma^2\bigr):=
\displaystyle\frac{\sqrt{2}}{{\sigma}},$}$\vspace*{2pt}
\cr
\sqrt{2}\sigma\beta, &\quad $\mbox{if $\beta\geq
\beta_c\bigl(\sigma^2\bigr).$}$}
\]
Corollary \ref{cor:main3} follows from the next result with the choice
$\sigma_1=\sigma_2$.

%
\begin{proposition}
\label{prop:perturbed-free-energy}
Let
$V_{12}:=\sigma_1^2\alpha+\sigma_2^2(1-\alpha)$. Then:
\begin{itemize}
\item\textit{Case} 1: If
$\sigma_1\leq\sigma_2$,
\[
\lim_{N\to\infty}f^{(\vec{\sigma},{\alpha})}_N(\beta)= f(\beta;
V_{12}).
\]
\item\textit{Case} 2:
If $\sigma_1\geq\sigma_2$,
\[
\lim_{N\to\infty}f_N^{(\vec{\sigma},{\alpha})}(\beta)= \alpha f
\bigl(\beta; \sigma_1^2\bigr) + (1-\alpha)f\bigl(\beta;
\sigma_2^2\bigr),
\]
\end{itemize}
where the convergence holds almost surely and in $L^1$.
\end{proposition}

The expressions are identical to the free energy of a GREM with two
levels. In case 1, it is reduced to a REM. The conditions can be
rewritten by defining a piecewise linear function of slopes $\sigma
_1^2$ and $\sigma_2^2$ on the intervals $[0,\alpha]$, $[\alpha,1]$,
respectively.
In case 1, this function fails to be concave. However, it is easily
verified that the effective parameters
define the concave hull of the function. The reader is referred to
\cite
{capocaccia-cassandro-picco} and \cite{bovier-kurkova1}
for more details on the concavity conditions which is very general for
the family of GREM models.
In case 1 there is one critical value for $\beta$, and in case 2 there
are two critical values for $\beta$ corresponding to the respective
$\beta_c(\sigma^2)$
of the two effective parameters $\sigma^2$. In case 1, the critical
$\beta$ is $\sqrt{2/V_{12}}$, whereas the two critical $\beta$'s are
$\sqrt{2}/\sigma_1$ and $\sqrt{2}/\sigma_2$
in case 2.


%
\begin{figure}[b]

\includegraphics{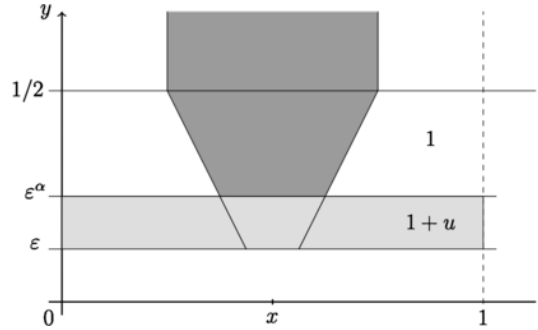}

\caption{The perturbed model where the variance parameter is $(1+u)$ on
the strip $[\varepsilon, \varepsilon^\alpha]$ where $\varepsilon=1/N$.}
\label{fig:perturbed-daviaud2}
\end{figure}


\subsection{The Bovier--Kurkova technique}\label{sec2.2}
The proof of Theorem \ref{thm:low-temperature} relies on determining
the overlap distribution of the original model from
the free energy of the perturbed ones. This approach has been used by
Bovier and Kurkova in the case of the GREM-type models \cite
{bovier-kurkova1,bovier-kurkova2}.

For $u \in(-1,1 )$ and $\alpha\in(0,1)$, consider the
field $(Y_x, x\in\mathcal X_\varepsilon)$ defined in \eqref{eqn:
perturbed Y} with the choice of parameters
$\vec{\sigma}=(1, (1+u))$; see Figure~\ref{fig:perturbed-daviaud2}.
(Recall that, for the sake of clarity, it is assumed that $N^\alpha$
and $N^{1-\alpha}$ are integers.)
The original Gaussian field $(X_x)$ is recovered at $u=0$.
Note that if $u>0$, the parameters correspond to the first case of
Proposition \ref{prop:perturbed-free-energy} and if $u<0$, to the second.
The field $Y$ can also be represented as follows:
%
%
\begin{equation}
\label{eqn: Y rep} Y_x= X_x + u \bigl(X_x-X_x(
\alpha) \bigr), \qquad 1 \le i \le N.
\end{equation}
The proof of the next lemma is a simple integration and is postponed to
the \hyperref[app]{Appendix}; see Appendix~\ref{proof-covariance-acc-perturbed}.
%
%
\begin{lemma}
\label{lem:covariance-acc-perturbed}
Fix $0<\varepsilon=1/N<1/2$, and $\alpha\in(0,1)$. Let
$\tilde X_x:=X_x-X_x(\alpha)$. Then, for $x\in\mathcal X_\varepsilon$
\[
\E\bigl[\tilde X_x^2\bigr] = \E[\tilde
X_xX_x]= (1 -\alpha) \log N, \qquad x\in\mathcal
{X}_\varepsilon,
\]
and, for $x,x'\in\mathcal X_\varepsilon$,
%
%
\begin{equation}
\label{eqn: x tilde covariance} \qquad\E[\tilde X_x X_{x'}]=\cases{
\bigl(q\bigl(x,x'\bigr)-\alpha
\bigr) \log N + O_N(1), & \quad $\mbox{if } \alpha< q\bigl(x,x'
\bigr) \le1,$
\vspace*{2pt}\cr
0, & \quad $\mbox{if } 0 \le q\bigl(x,x'\bigr) \le\alpha,$}
\end{equation}
where $O_N(1)$ is a term uniformly bounded in $N$, and we recall that
$\Vert x -x' \Vert= \varepsilon^{q(x,x')}$.
\end{lemma}

This result and a Gaussian integration by parts yield an important lemma.
%

\begin{lemma}
\label{lem:perturbed-Ipp} For all $\alpha\in(0,1)$, we have
\[
\beta\int_{\alpha}^{1} x_\beta^{(N)}(s)
\,\d s + o_N(1)=\frac{1}{\log
N}{\mathbb E}\biggl[\sum
_{x\in\mathcal X_\varepsilon} G_{\beta, N}(x) \bigl( X_x-
X_x(\alpha)\bigr) \biggr],
\]
where $o_N(1)$ stands for a term that goes to $0$ as $N$ goes to
$\infty$.
\end{lemma}

\begin{pf}
Fix $\varepsilon=1/N$ and $\alpha\in(0,1)$. Note that $(\tilde X_x;
(X_{x'}, x'\in\mathcal X_\varepsilon))$ is a Gaussian vector of $N+1$
variables.
Therefore, Gaussian integration by parts (see Lemma \ref{lem:
GaussianIbP}) yields, for all
$x\in\mathcal X_\varepsilon$,
\begin{eqnarray*}
\beta^{-1} {\mathbb E}\biggl[ \frac{\tilde X_x\ee^{\beta X_x}}{\sum
_{x'\in
\mathcal X_\varepsilon} \ee^{\beta X_{x'}}} \biggr] &=&- \sum
_{ x'\in\mathcal{X}_\varepsilon} {\mathbb E}[\tilde X_x X_{x'} ]
{\mathbb E}
\biggl[\frac{\ee^{\beta(X_x+X_{x'})}}{ (\sum_{z\in
\mathcal{X}_\varepsilon} \ee^{\beta X_z} )^2} \biggr]
\\
&&{} + {\mathbb E}[\tilde X_xX_x ] {\mathbb E}\biggl[
\frac{\ee^{\beta X_x}}{\sum_{z\in\mathcal{X}_\varepsilon} \ee
^{\beta X_z}} \biggr].
\end{eqnarray*}
Lemma \ref{lem:covariance-acc-perturbed} and elementary manipulations imply
\begin{eqnarray*}
&& (\beta\log N)^{-1} {\mathbb E}\biggl[ \sum_{x\in\mathcal
X_\varepsilon}
\tilde X_xG_{\beta, N}(x) \biggr]
\\
&&\qquad= \sum_{x,x'\in\mathcal X_\varepsilon} \biggl( \int_{\alpha}^{1}
\mathbf{1}_{\{q(x,x')\le s\}} \,\d s \biggr) {\mathbb E}\bigl[
G_{\beta, N}(x)
G_{\beta, N}\bigl(x'\bigr) \bigr] + O \biggl(\frac{1}{\log N}
\biggr)
\\
&&\qquad= \int_{\alpha}^{1} {\mathbb E}\bigl[ G_{\beta, N}^{\times2}
\{q_{12} \le s \} \bigr] \,\d s + O \biggl(\frac{1}{\log N} \biggr),
\end{eqnarray*}
which concludes the proof of the lemma.
\end{pf}

\begin{pf*}{Proof of Theorem \ref{thm:low-temperature}}
Fix $\beta>\beta_c=\sqrt{2}$.
Write $Z_N^{(u, \alpha)}(\beta)$ for the partition function \eqref
{eq:partition-function-perturbed} for the choice $\vec{\sigma}= (1, (1+u))$.
Direct differentiation and equation~\eqref{eqn: Y rep} give
\[
\frac{\d}{\d u} \bigl( {\mathbb E}\log Z_N^{(u, \alpha)}(\beta)
\bigr)_{u=0} = \beta{\mathbb E}\biggl[\sum_{x\in\mathcal
X_\varepsilon}
\bigl( X_x- X_x(\alpha)\bigr)G_{\beta, N}(x)
\biggr],
\]
which, together with Lemma \ref{lem:perturbed-Ipp}, yields
%
%
\begin{equation}
\label{eqn: convex 1} \int_{\alpha}^{1} x_\beta^{(N)}(s)
\,\d s = \beta^{-2} (\log N)^{-1} \frac{\d}{\d u} \bigl( {\mathbb E}
\log Z_N^{(u, \alpha)}(\beta) \bigr)_{u=0}
+o_N(1).
\end{equation}
Observe that ${\mathbb E}f_N^{(u, \alpha)}(\beta) = (\log N)^{-1}
{\mathbb E}\log
Z_N^{(u, \alpha)}(\beta)$ is a convex function of $u$.
Moreover, by Proposition \ref{prop:perturbed-free-energy}, ${\mathbb E}f_N^{(u,
\alpha)}(\beta)$ converges.
The limit, that we denote $f^{(u, \alpha)}(\beta)$, is also convex in
the parameter $u$.
In particular, by a standard result of convexity (see, e.g.,
Proposition I.3.2 in \cite{simon}), at every point of
differentiability, the derivative of the limit equals the limit
of the derivative
%
%
\begin{eqnarray}
\label{eqn: convex 2} \lim_{N \to\infty} \frac{\d}{\d u} {\mathbb E}
f_N^{(u, \alpha)}(\beta) = \frac{\d}{\d u} f^{(u, \alpha)}(
\beta)
\nonumber
\\[-8pt]
\\[-8pt]
 \eqntext{\mbox{$\forall u$ where $u\mapsto f^{(u, \alpha)}(\beta)$ is
differentiable}.}
\end{eqnarray}
We show $f^{(u, \alpha)}(\beta)$ is differentiable at $u=0$.
The derivative can be computed by Proposition \ref{prop:perturbed-free-energy}.
For $u$ small enough, $\beta$ is larger than all critical $\beta$'s. Thus
%
%
\begin{equation}
\frac{\d}{\d u} f^{(u, \alpha)}(\beta)=\cases{ %
\displaystyle\sqrt{2} \beta\frac{(1-\alpha)(1+u)}{\sqrt{\alpha+
(1-\alpha)(1+u)^2}}, & \quad$\mbox{if } u>0,$
\vspace*{2pt}\cr
\sqrt{2} \beta(1-\alpha),&\quad $\mbox{if } u<0.$}
\end{equation}
From this, it is easily verified that $f^{(u, \alpha)}(\beta)$ is
differentiable at $u=0$ and
%
%
\begin{equation}
\label{eqn: convex 3} \frac{\d}{\d u} \bigl(f^{(u, \alpha)}(\beta)
\bigr)_{u=0}=\sqrt{2} \beta(1- \alpha).
\end{equation}
Equations \eqref{eqn: convex 1}, \eqref{eqn: convex 2} and \eqref{eqn:
convex 3} together imply
%
%
\begin{equation}
\label{eqn: weak conv} \lim_{N\to\infty} \int_{\alpha}^{1}
x^{(N)}_{\beta}(s) \,\d s = \frac
{\sqrt{2}}{\beta} (1-\alpha)\qquad \mbox{for all $\alpha\in(0,1)$}.
\end{equation}
Therefore, any weak limit $x_\beta$ must satisfy $x_\beta(\alpha)
\le
\sqrt{2}/\beta$ for any point of continuity $\alpha<1$,
since $x_\beta$ is nondecreasing.
If there exists $0<\alpha<1$ such that $x_\beta(\alpha)<\frac{\sqrt
{2}}{\beta}$, there would be a contradiction with \eqref{eqn: weak
conv}, since by right-continuity and monotonicity of $x_\beta$ we could
find $\alpha'>\alpha$
such that
\[
\lim_{N\to\infty} \int_{\alpha}^{\alpha'}
x^{(N)}_{\beta}(s) \,\d s < \frac{\sqrt{2}}{\beta} \bigl(
\alpha'-\alpha\bigr).
\]
This proves that any weak limit $x_\beta$ of $(x^{(N)}_{\beta}, N\in
\N
)$ is the same and equals $\frac{\sqrt{2}}{\beta}$ on $(0,1)$.
The subsequential limits being the same, this proves in particular
convergence of the sequence to the desired distribution function.
\end{pf*}

\subsection{A spin-glass approach to Poisson--Dirichlet variables}
\label{sect: SG->PD}
In this section, the link between Theorems \ref{thm:low-temperature}
and~\ref{thm:main4} is explained.
The technique, inspired from the study of spin glasses in particular
\cite{arguin-chatterjee}, is general and is of independent interest to
prove convergence to Poisson--Dirichlet statistics.

The first step is to find a good space for the convergence of $G_{\beta,N}$.
Let $C$ be the compact metric space of $\N\times\N$ covariance matrices
with $1$ on the diagonal endowed with the product topology on the entries.
For a given $N$, consider the mapping
\begin{eqnarray*}
\mathcal X_N^{\times\infty} &\to& C,
\\
(x_l, l\in\N) &\mapsto& R^{(N)},
\end{eqnarray*}
where for $l,l'\in\N$
\[
R^{(N)}_{l,l'}:= \cases{q_{ll'}=q
\bigl(x_l,x_l'\bigr), &\quad $\mbox{if $l\neq
l'$}$\vspace*{2pt}
\cr
1, &\quad $\mbox{if $l=l'$.}$}
\]
Consider the probability measure $\E G_{\beta,N}^{\times\infty}$ on
$\mathcal X_N^{\times\infty}$.
The push-forward of this probability measure under the above mapping
defines a random element of $C$ that we denote $\vec{R}^{(N)}$.
Since each point is sampled independently from the same measure, the
law of $\vec{R}^{(N)}$ is \textit{weakly exchangeable}, that is, for
any permutation $\pi$ of a finite number of indices,
\[
\bigl(\vec{R}^{(N)}_{\pi(l)\pi(l')}\bigr) \,\equilaw\,\bigl(
\vec{R}^{(N)}_{ll'}\bigr).
\]

The sequence of random matrices $(\vec{R}^{(N)}, N\in\N)$ is tight by
Prokhorov's theorem since the space $C$ is a compact metric space.
Hence, there exists a subsequence $\{\vec{R}^{(N_m)}\}_{m\in\N}$ that
converges weakly. Denote the subsequential limit by $\vec{R}$.
Observe that $\vec{R}$ is also weakly exchangeable since the mappings
on $C$ induced by a finite permutation is continuous.
Therefore, by the\vspace*{1pt}  representation theorem of Dovbysh and Sudakov \cite
{dovbysh-sudakov}, $\vec{R}$ is constructed like $\vec{R}^{(N)}$ by
sampling from a random measure.
Precisely, the theorem states that there exists a random
probability\vadjust{\goodbreak}
measure $\mu_\beta$ on a Hilbert space $\mathcal H$, with law $P$ and
corresponding expectation $E$,
such that the random matrix $\vec{R}$ has the same law as the Gram
matrix of a sequence of vectors $(v_l,l\in\N)$ that are sampled under
$E\mu_\beta^{\times\infty}$.
[In other words, the vectors $(v_l,l\in\N)$ are i.i.d. conditionally on
$\mu_\beta$.]
The equality in law can be expressed as follows: for any continuous
function $F$ on $C$,
%
%
\begin{equation}
\label{eqn: mu} \lim_{m\to\infty} \E G_{\beta,N_m}^{\times\infty}
\bigl[F(q_{ll'}) \bigr]= E\mu_{\beta}^{\times\infty}\bigl[F
\bigl(v_l \cdot v_l'\bigr)\bigr].
\end{equation}
Note that, since $ q(x,x')\leq1$, the random measure $\mu_{\beta}$ is
supported on the unit ball.
The first consequence of Theorem \ref{thm:low-temperature} is that for
any subsequential limit $\mu_\beta$,
%
%
\begin{eqnarray}
\label{eqn: x mu} E \bigl[ \mu_{\beta}^{\times2} \{v_1
\cdot v_{2} \le q \} \bigr] &=&\lim_{N\to\infty} {\mathbb E}\bigl[
G_{\beta, N}^{\times2} \{ q_{12} \le q \} \bigr]
\nonumber
\\[-8pt]
\\[-8pt]
\nonumber
& =&
\frac{\beta_c}{\beta} 1_{[0,1)}(q)+ 1_{\{
1\}}(q).
\end{eqnarray}
The first equality is obtained by bounding $1_{[0,q]}(q_{ll'})$ by
continuous functions on $q_{ll'}$ above and below and by applying
\eqref
{eqn: mu}.
In view of equations \eqref{eqn: mu} and \eqref{eqn: x mu}, we see the
random measures $\mu_\beta$ as limit points of $(G_{\beta,N})_{N\in
\N}$.

The main ingredient to prove Poisson--Dirichlet statistics is a general
property of the Gibbs measure $(G_{\beta,N}(x), x \in\mathcal X_N)$
of centered Gaussian fields known as the Ghirlanda--Guerra identities.
They were introduced in \cite{ghirlanda-guerra} and were proved in a
general setting by Panchenko \cite{panchenko_gg}.
%
%
\begin{theorem}
\label{thm: GG}
Let $\mu_\beta$ be a subsequential limit of $(G_{\beta,N})_{N\in\N
}$  in\break
the sense of~\eqref{eqn: mu}.
Then for any $s\in\N$ and any continuous functions
$F\dvtx\break [-1,1]^{s(s-1)/2}\to\R$
%
%
\begin{eqnarray}
\label{eqn: pre GG} 
E \mu_\beta^{\times s+1}
\bigl[v_1\cdot v_{s+1} F(v_l\cdot
v_{l'}) \bigr]
&=&\frac{1}{s} E \mu_\beta^{\times2} [v_1\cdot
v_2 ] E \mu_\beta^{\times s} \bigl[F(v_l
\cdot v_{l'}) \bigr]
\nonumber
\\[-8pt]
\\[-8pt]
\nonumber
&&{} + \frac{1}{s}\sum
_{k=2}^s E \mu_\beta^{\times s}
\bigl[v_1\cdot v_k F(v_l\cdot
v_{l'}) \bigr]. \nonumber 
\end{eqnarray}
\end{theorem}

\begin{pf}
Recall that we write $G_{\beta,N}^{\times s}$ for the product measure
on $\mathcal X_N^s$.
Also for $(x_1,\ldots, x_s)\in\mathcal X_N^s$, the overlaps
$q(x_l,x_l')$, $1\leq l,l'\leq s$, are denoted $q_{ll'}$.
In a similar way, we write $X_1$ for the field $X_{x_1}$ of the first
point sampled from $G_{\beta,N}$.
It is shown in \cite{panchenko_gg} that, for any $\beta$ where the free
energy $f(\beta)$ is differentiable, the following concentration holds:
%
%
\begin{equation}
\label{eqn: sa} \lim_{N\to\infty}\frac{1}{\log N} \E
G_{\beta,N} \bigl[\bigl | X_{1} - \E G_{\beta,N}(X_{1})
\bigr| \bigr]=0.
\end{equation}
Note that by Corollary \ref{cor:main3}, differentiability holds at all
$\beta$ for the Gaussian field considered.
Since the function $F$ is bounded, \eqref{eqn: sa} implies
%
%
\begin{equation}
\label{eq:GG-double+}\quad  \lim_N \frac{1}{\log N} \bigl( \E
G_{\beta,N}^{\times s} \bigl[X_1 F(q_{ll'})
\bigr]- \E G_{\beta,N}[X_1] \E G_{\beta,N}^{\times
s}
\bigl[F(q_{ll'}) \bigr] \bigr)=0.\vadjust{\goodbreak}
\end{equation}
The two terms can be evaluated by Gaussian integrations by part (see
Lemma \ref{lem: GaussianIbP}),
%
%
\begin{equation}
\label{eq-GG-Ipp1+} \frac{1}{\beta\log N} \E G_{\beta,N}[X_1]= 1-
\E
G_{\beta,N}^{\times
2}[q_{12}] + O \biggl(\frac{1}{\log N}
\biggr)
\end{equation}
and
%
%
\begin{eqnarray}
\label{eq-GG-Ipp2+} && \frac{1}{\beta\log N} \E G_{\beta,N}^{\times s}
\bigl[X_1 F(q_{ll'})\bigr]
\nonumber
\\
&&\qquad =-s \E G_{\beta,N}^{\times s+1} \bigl[q_{1,s+1}
F(q_{ll'})\bigr]+ \sum_{1 \le k
\le s} \E
G_{\beta,N}^{\times s} \bigl[q_{1k} F(q_{ll'})
\bigr] \\
&&\qquad\quad{}+ O \biggl(\frac
{1}{\log N} \biggr).\nonumber
\end{eqnarray}
Finally recalling \eqref{eq:GG-double+} and assembling \eqref
{eq-GG-Ipp1+}--\eqref{eq-GG-Ipp2+} yields the Ghirlanda--Guerra
identities (see equation (16) in \cite{ghirlanda-guerra}),
%
%
\begin{eqnarray}
\label{eqn: GG proof} %
&&\E G_{\beta,N}^{\times s+1}
\bigl[q_{1,s+1} F(q_{ll'}) \bigr]
\nonumber
\\
&&\qquad=\frac{1}{s} \E G_{\beta,N}^{\times2} [q_{12} ] \E
G_{\beta,N}^{\times s} \bigl[F(q_{ll'}) \bigr] +
\frac{1}{s}\sum_{k=2}^s\E
G_{\beta,N}^{\times s} \bigl[q_{1k} F(q_{ll'})
\bigr]\\
&&\qquad\quad{} +o_N(1).\nonumber 
\end{eqnarray}
[Note that the term for $k=1$ cancels with the 1 since $q_{11}=1+o_N(1)$.]
In particular, for any subsequential limit $\mu_\beta$ of $(G_{\beta
,N})_N$ in the sense of \eqref{eqn: mu}, one obtains~\eqref{eqn: pre
GG} by taking the limit $N\to\infty$ and applying the definition of
convergence in the sense of \eqref{eqn: mu}.
\end{pf}

Equation \eqref{eqn: x mu} and the Ghirlanda--Guerra identities imply
that $\mu_\beta$ is atomic.
%
%
\begin{corollary}
\label{cor: atomic}
Let $\mu_\beta$ be a subsequential limit of $(G_{\beta,N})_{N\in\N
}$ in
the sense of \eqref{eqn: mu}.
Then there exist random weights $\xi=(\xi_i;i\in\N)_\downarrow$ with
$\xi_i\geq0$, $\sum_{i\in\N}\xi_i=1$
and orthonormal vectors $(e_i; i\in\N)\subset\mathcal H$
such that
\[
\mu_\beta=\sum_{i\in\N}\xi_i
\delta_{e_i},\qquad P\mbox{-a.s.}
\]
Moreover, from \eqref{eqn: x mu}, $E[\sum_{i\in\N}\xi_i^2]=1-
\frac
{\beta_c}{\beta}$.
\end{corollary}

\begin{pf}
Let $(v_l,l\in\N)$ be a sequence sampled from $E\mu_\beta^{\times
\infty}$.
From $(v_l,l\in\N)$, we reconstruct $\mu_\beta$ up to isometry.
For a fixed $l$ consider the sequence $(v_l\cdot v_{l'}, l' > l)$. This
is a sequence of 0's and 1's by \eqref{eqn: x mu}.
We first show that, almost surely, for every $l\in\N$, there exists
$l'>l$ such that $v_l\cdot v_{l'}=1$;
in particular, since all vectors are in the\vadjust{\goodbreak} unit ball, $v_l=v_{l'}$ and
$\|v_l\|=1$.
For this, we proceed as in Lemma~1 in \cite{panchenko_gg2}.
Write $F_s(v_{l}\cdot v_{l'})=\prod_{l=2}^{s}(1-v_{1}\cdot v_{l})$. In
other words, $F_s(v_{l}\cdot v_{l'})$ is $1$
if $v_{1}\cdot v_{l}=0$ for $l=2,\ldots, s$, otherwise it is $0$.
Denote for short $\alpha=1-E\mu_\beta^{\times2}[v_1\cdot v_2]$.
Equation \eqref{eqn: pre GG} implies
\begin{eqnarray*}
&&E\mu_\beta^{\times s+1}\{v_1
\cdot v_l=0, 2\leq l\leq s+1\}
\\
&&\qquad= E\mu_\beta^{\times s+1}\bigl[(1-v_{1}\cdot
v_{s+1})F_s(v_{l}\cdot v_{l'})\bigr]
\\
&&\qquad= \frac{\alpha}{s}E\mu_\beta^{\times s}\{v_1\cdot
v_l=0, 2\leq l\leq s\} +\frac{1}{s}\sum
_{l=2}^sE\mu_\beta^{\times s}
\{v_1\cdot v_l=0, 2\leq l\leq s\}
\\
&&\qquad=\frac{s-1+\alpha}{s}E\mu_\beta^{\times s}\{v_1\cdot
v_l=0, 2\leq l\leq s\}= \frac{(s-1+\alpha)\cdots(1+\alpha) \alpha
}{s!}, 
\end{eqnarray*}
where the last equality is obtained by induction. The last term goes to
$0$ as $s\to\infty$ since $\alpha<1$, hence
\[
E\mu_\beta^{\times\infty}\{v_1\cdot v_l=0, l
\geq2\}=0,
\]
from which we deduce that, $P$-a.s., $\mu_\beta^{\times\infty}\{
v_1\cdot v_l=0, l\geq2\}=0$ and then that, for $\mu_\beta$-almost
all $v$,
\[
\mu_\beta^{\times\infty}\{v\cdot v_l=0, l\geq2\} =0.
\]
Since the vectors $v_l$ are i.i.d. $\mu_\beta$-sampled, it follows
that, $P$-a.s., for $\mu_\beta$-almost all $v$, $\mu_\beta(v\cdot
v_1=0)<1$, thus $\mu_\beta(v\cdot v_1=1)>0$
as claimed.

By the reasoning above, a vector that is sampled once in $(v_l,l\in\N)$
is sampled infinitely many times $E\mu_\beta^{\times\infty}$-a.s.
Moreover, since the vectors are conditionally i.i.d., for $l\in\N$, the
following limit exists and must be nonzero:
%
%
\begin{equation}
\label{eqn: weight} \lim_{n\to\infty}\frac{1}{n}\sum
_{j=l+1}^{l+n} 1_{v_l}(v_j)>0,\qquad
\mbox{$E\mu_\beta^{\times\infty}$-a.s.}
\end{equation}
In particular, every sampled vector $v_l$ is an atom a.s. and its
weight is measurable with respect to $(v_l,l\in\N)$.
Moreover, if $v_l\neq v_{l'}$, then $v_l\cdot v_{l'}=0$ $E\mu_\beta
^{\times\infty}$-a.s. Therefore the atoms are orthogonal.
It remains to consider the different atoms without repetitions and
reorder the weights.
Let $e_1=v_{1}$, $e_2=v_{l_2}$ where $l_2=\inf\{l\geq1\dvtx v_l\cdot
e_1=0\}$, $e_3=v_{l_3}$ where $l_3=\inf\{l\geq l_2\dvtx v_l\cdot e_i=0,
i=1,2\}$
, and so forth.
By construction, $(e_j, j\geq1)$ are orthonormal vectors. (The
collection is not necessarily infinite at this point.)
We can assign to each vector $e_j$ its weight $\mu_\beta(\{e_j\})$ by
\eqref{eqn: weight}.
The collection can then be ordered in decreasing order to get the result.

The fact that $E[\sum_{i\in\N}\xi_i^2]=1- \frac{\beta_c}{\beta}$ is
straightforward from \eqref{eqn: x mu}.
\end{pf}

To finish the proof of Theorem \ref{thm:main4}, it remains to show that
the random weights $\xi$ are distributed like a Poisson--Dirichlet
variable of parameter $\frac{\beta_c}{\beta}$.
In fact, the parameter is already determined by Corollary \ref{cor:
atomic}, since for a Poisson--Dirichlet variable $\xi'$ of parameter
$x$, $E[\sum_k (\xi_k')^2]=1-x$ holds; see, for example, Corollary~2.2
in \cite{ruelle}.
This will also imply that for any converging sequence of $(G_{\beta
,N})$ in the sense of \eqref{eqn: mu}, the limit is the same. In
particular, it implies convergence of the whole sequence by compactness.

To prove the Poisson--Dirichlet statistics of the weights $\xi$, we use
the following characterization theorem of the law; see \cite
{talagrand}, page 22 for details.
Define for all $m\in\N$ the joint moments of the weights
%
%
\begin{equation}
\label{eqn: S} S(n_1,\ldots,n_m)=E \sum
_{k_1,\ldots,k_m} \xi_{k_1}^{n_1}\cdots\xi
_{k_m}^{n_m}\qquad \mbox{for $n_1,\ldots,n_m
\ge1$}.
\end{equation}
The collection of $S(n_1,\ldots,n_m)$, $m\in\N$, determines the law of a
random mass-partition, that is, a random variable on ordered sequences
$1\geq r_1 \geq r_2 \geq\cdots\geq0$ with $\sum_{i\in\N}r_i\leq1$.
If $\xi$ is a Poisson--Dirichlet variable, it is shown in \cite
{talagrand}, Proposition 1.2.8, that the moments satisfy the recursion relations
\begin{eqnarray}\qquad
\label{eqn: GG} %
S(n_1+1,\ldots,n_m)&=&
\frac{S(2)}{s} S(n_1,\ldots,n_m)+\frac{n_1-1}{s}
S(n_1,\ldots,n_m)
\nonumber
\\[-8pt]
\\[-8pt]
\nonumber
 &&{} + \sum
_{2\leq l \leq m} \frac{n_l}{s}S(n_1+n_l,n_2,
\ldots, n_{l-1}, n_{l+1},\ldots, n_m),
\end{eqnarray}
where $s=n_1+\cdots+n_m$. It is not hard to verify that all moments
$S(n_1,\ldots,n_m)$ (and thus the law of $\xi$) are determined by
recursion from $S(2)$ and the identities~\eqref{eqn: GG}.

It turns out that these identities are satisfied by $\xi$ defined by
Theorem \ref{thm: GG} and Corollary \ref{cor: atomic}.
%
%
\begin{theorem}
\label{thm: pd}
Let $\xi$ be a random mass-partition satisfying the assumptions of
Corollary \ref{cor: atomic}.
The moments $S(n_1,\ldots,n_m)$ of $\xi$ satisfy \eqref{eqn: GG} for any
$m\in\N$ and any $n_1,\ldots, n_m\in\N$.
In particular, $\xi$ has the law of a Poisson--Dirichlet variable of
parameter $1-S(2)$.
\end{theorem}
\begin{pf}
To deduce \eqref{eqn: GG} from \eqref{eqn: pre GG}, we follow \cite
{talagrand}, pages 24--25.
The set $\{1,\ldots,s\}$ can be decomposed into the disjoint union of
sets $I_1, \ldots, I_m$ with $\vert I_j \vert= n_j$ for all $1 \le j
\le m$.
Consider the functions $(F_j)_{1 \le j \le m}$ given by $F_j(\delta
_{k_lk_{l'}}):=\prod_{k_l,k_{l'} \in I_j} \delta_{k_lk_{l'}}$ and
define\vspace*{2pt} $F:=\prod_{1 \le j \le m} F_j$.
Then elementary manipulations imply \eqref{eqn: GG}.
Note that the second term on the right-hand side of \eqref{eqn: pre GG}
yields the last two terms of \eqref{eqn: GG}.
\end{pf}

\section{High points of the perturbed models}
\label{sec:perturbed-maximum}
In this section, the log-number of high points at a given level is
computed for the perturbed models introduced in Section~\ref{sect:outline}.
The focus is on the Gaussian field introduced in Section~\ref{sect:
perturbed}, though the technique
applies to any perturbed model with a finite number of parameters.
The free energies of the models are computed in Section~\ref{sect:
free energy}.

Let $Y=(Y_x, x\in\mathcal X_\varepsilon)$ be the Gaussian field
introduced in Section~\ref{sect: perturbed}.
Recall the notation and the two choices of parameters in Proposition
\ref{prop:perturbed-free-energy}:
%
%
\begin{eqnarray}
&&\mathit{Case}\ 1\dvtx\qquad  \sigma_1\leq
\sigma_2;
\nonumber
\\[-8pt]
\\[-8pt]
\nonumber
&&\mathit{Case}\ 2\dvtx\qquad \sigma_1\geq\sigma_2. 
\end{eqnarray}
Define also as before $V_{12}:=\sigma_1^2\alpha+\sigma_2^2(1-\alpha)$.
%
%
\begin{proposition}
\label{prop:perturbed-maximum}
\[
\lim_{N \to\infty} \p\Bigl( \max_{x\in\mathcal X_\varepsilon}
Y_x \ge\sqrt{2} \gamma_{\mathrm{max}} \log N \Bigr) =0,
\]
where
\[
\gamma_{\mathrm{max}}=\gamma_{\mathrm{max}}(\vec{\sigma},\alpha):= \cases{
\sqrt{V_{12}}, &\quad $\mbox{for case 1;}$ \vspace*{2pt}
\cr
\sigma_1\alpha+\sigma_2(1-\alpha), &\quad $\mbox{for case 2.}$}
\]
\end{proposition}

%
\begin{proposition}
\label{prop:perturbed-highpoints}
Let $\mathcal{H}_N^{Y}(\gamma):= \{ x\in\mathcal X_\varepsilon:
Y_x \ge\sqrt{2} \gamma\log N \}$ be the set of $\gamma$-high points.
Then, for all $0< \gamma< \gamma_{\mathrm{max}}$,
\[
\lim_{N \to\infty} \frac{\log\vert\mathcal{H}_N^{Y}(\gamma
)\vert
}{\log N}=\mathcal{E}^{(\vec{\sigma},\alpha)}(
\gamma) \qquad\mbox{in probability,}
\]
where in case 1,
\[
\mathcal{E}^{(\vec{\sigma},\alpha)}(\gamma):=1- \frac{\gamma^2}{V_{12}};
\]
and in case 2,
\[
\mathcal{E}^{(\vec{\sigma},\alpha)}(\gamma):= \cases{\displaystyle 1- \frac
{\gamma^2}{V_{12}}, &\quad $\mbox{if $ \displaystyle\gamma< \frac{V_{12}}{\sigma_1}$}$, \vspace*{2pt}
\cr
\displaystyle(1-\alpha) -
\frac{(\gamma-\sigma_1\alpha)^2}{\sigma_2^2(1-\alpha)}, &\quad $\mbox{if $ \displaystyle\gamma\geq\frac{V_{12}}{\sigma_1}$.}$}
\]

Moreover, for any $\mathcal{E}<\mathcal{E}^{(\vec{\sigma},\alpha
)}(\gamma)$, there exists $c$ such that
\[
\p\bigl(\bigl\vert\mathcal{H}_N^{Y}(\gamma) \bigr\vert\le
N^{\mathcal{E}} \bigr) \le\exp\bigl\{- c (\log N)^2\bigr\}.
\]
\end{proposition}

\subsection{Proof of Proposition \texorpdfstring{\protect\ref{prop:perturbed-maximum}}{3.1}}
The proof of case 1 is by a union bound,
\[
\p\Bigl( \max_{x\in\mathcal X_\varepsilon} Y_x \geq\sqrt{2}
\gamma
_{\mathrm{max}} \log N \Bigr)\leq N \p( Y_x \geq\sqrt{2}
\gamma_{\mathrm{max}} \log N ),
\]
which goes to zero by a Gaussian estimate; see Lemma \ref{lem: gaussian}.
For case 2, we construct a Gaussian field with hierarchical
correlations that dominates $Y$ at the level of the covariances.
The result will follow by comparison using Slepian's lemma.

Notice that if $\varepsilon< \|x-x'\| \leq\varepsilon^{\alpha}$,
the corresponding cone-like sets for $Y_x$ and~$Y_{x'}$ in $ \mathcal
C^+$ intersect between the lines $y=\varepsilon$ and $y=\varepsilon
^{\alpha}$.
Therefore the covariance of the variables satisfies, writing $\ell
:=\Vert x - x' \Vert$,
\begin{eqnarray*}
\E[Y_x Y_{x'}] &=& \sigma_2^2
\int_{\ell}^{\varepsilon^{\alpha}} \frac
{y-\ell}{y^2} \,\d y +
\sigma_1^2 \biggl( \int_{\varepsilon^{\alpha
}}^{1/2}
\frac{y-\ell}{y^2} \,\d y + \int_{1/2}^{\infty}
\frac{1/2-\ell
}{y^2} \,\d y \biggr)
\\
&\geq& \sigma_1^2 \biggl( \log\frac{1/2}{\varepsilon^{\alpha}} -1
\biggr).
\end{eqnarray*}
By applying the same reasoning when $ \varepsilon^{\alpha} < \|x-x'\|
\le1/2$, one obtains the following lower bound for the covariance:
%
%
\begin{equation}
\label{eqn: cov bound} \E[Y_x Y_{x'}]\geq\cases{0, &\quad
$\mbox{if $\bigl\|x-x'\bigr\|> \varepsilon^{\alpha},$}$\vspace*{2pt}
\cr
\displaystyle\sigma_1^2 \biggl(\log\frac{1/2}{\varepsilon^{\alpha}}-1 \biggr
), &\quad  $\mbox{if $\varepsilon< \bigl\|x-x'\bigr\| \leq\varepsilon^{\alpha}.$}$}
\end{equation}
Equation \eqref{eqn: cov bound} is used to construct a Gaussian field
$\tilde Y$.
Define the map
\begin{eqnarray*}
\pi\dvtx\mathcal X_{\varepsilon} &\to&\mathcal X_{\varepsilon
^{\alpha}},
\\
x &\mapsto&\pi(x),
\end{eqnarray*}
where $\pi(x)$ is the unique $y\in\mathcal X_{\varepsilon^{\alpha}}$
such that $\| x-y \|\leq\frac{\varepsilon^{\alpha}}{2}$.
(If $\| x-y \|= \frac{\varepsilon^{\alpha}}{2}$, there are two
possibilities for $y$. We take the right point.)
The pre-image\vspace*{1.5pt} of $y \in\mathcal X_{\varepsilon^{\alpha}}$ under $\pi$
are exactly the points in $ \mathcal X_{\varepsilon}$
that are at a distance less than $ \frac{\varepsilon^{\alpha}}{2}$ from
$y$. One can think of $\pi(x)$ as the \textit{ancestor} of $x$ at the
scale $\varepsilon^{\alpha}$.

Consider the following Gaussian variables
%
%
\begin{eqnarray}
&&\bigl(g^{(1)}_{x}, x\in\mathcal
X_{\varepsilon^{\alpha}}\bigr) \qquad\mbox{i.i.d. Gaussians of variance $
\sigma_1^2\alpha\log N -\sigma_1^2
\log2-\sigma_1^2$,}
\nonumber\hspace*{-35pt}
\\[-8pt]
\\[-8pt]
\nonumber
&&\bigl(g^{(2)}_{x}, x\in\mathcal X_{\varepsilon}\bigr) \qquad
\mbox{i.i.d. Gaussians of variance $\sigma_2^2(1-\alpha)
\log N+2\sigma_1^2$.}\hspace*{-35pt}
\end{eqnarray}
These two families are also assumed independent.
Then, the field $\tilde Y$ is defined, using the map $\pi$ above and
the Gaussian random variables $g_x^{(i)}$, by
%
%
\begin{equation}
\label{eqn: g} \tilde Y_x= g_{\pi(x)}^{(1)} +
g_x^{(2)}.
\end{equation}

This construction and equation \eqref{eqn: cov bound} directly imply
the following comparison lemma.
%
%
\begin{lemma}
\label{lem: comparison}
%
%
\begin{eqnarray}
\E\bigl[\tilde Y_x^2\bigr]&=&\E
\bigl[Y_x^2\bigr]\qquad \mbox{$\forall x\in\mathcal
X_\varepsilon$},
\nonumber
\\[-8pt]
\\[-8pt]
\nonumber
\E[\tilde Y_x \tilde Y_y]&\leq&\E[Y_x
Y_y]\qquad \mbox{$\forall$ $x\neq y$, $x,y\in\mathcal
X_\varepsilon$.}
\end{eqnarray}
\end{lemma}

The following corollary is a straightforward consequence of the above
lemma and Slepian's lemma; see Corollary 3.12 in \cite{ledoux-talagrand}.
%
%
\begin{corollary}
\label{cor: bound max}
For any $\lambda>0$,
%
%
\begin{equation}
\p\Bigl( \max_{x\in\mathcal X_\varepsilon} Y_x \ge\lambda\Bigr)
\leq\p\Bigl( \max_{x\in\mathcal X_\varepsilon} \tilde Y_x \ge
\lambda
\Bigr).
\end{equation}
\end{corollary}

The Gaussian field $\tilde Y$ is almost identical to a GREM model with
two levels with parameters $0<\alpha<1$ and $\sigma_1,\sigma_2$; see,
for example, \cite{derrida,bovier-kurkova1}.
In fact the only aspect different from an exact GREM are the terms of
order one in the variances of the Gaussian random variables $g_x^{(i)}$'s.
However, these do not affect the first order of the maximum.
The proof of Proposition \ref{prop:perturbed-maximum} is concluded by
the following standard GREM result.
The proof of the lemma is not hard and is omitted for conciseness.
The reader is referred to Theorem 1.1 in \cite{bovier-kurkova1} where a
stronger result on the maximum is given and to \cite
{bolthausen-sznitman}, Lecture 9, for more details on the free energy
and on the log-number of high points of a two-level GREM.

%
\begin{lemma}
\label{lem:grem}
Let $\tilde Y$ be the Gaussian field constructed above. Then
\[
\p\Bigl( \max_{x\in\mathcal X_\varepsilon} \tilde Y_x \ge\sqrt{2}
\gamma_{\mathrm{max}}\log N \Bigr) \to0,\qquad N\to\infty,
\]
where $\gamma_{\mathrm{max}}$ is defined in Proposition \ref{prop:perturbed-maximum}.
\end{lemma}

\subsubsection{Proof of the upper bound in Proposition \texorpdfstring{\protect\ref{prop:perturbed-highpoints}}{3.2}}
The goal is to get an upper bound in probability) for $\vert\mathcal
{H}^Y_N(\gamma) \vert$ where $\mathcal{H}^Y_N(\gamma)=\{x \in
\mathcal
X_\varepsilon\dvtx Y_x\geq\sqrt{2}\gamma\log N\}$.

In case 1, a first moment computation gives the result. Indeed, a
Gaussian estimate (see Lemma \ref{lem: gaussian}) gives
\[
\E\bigl[\bigl\vert\mathcal{H}_N^{ Y}(\gamma) \bigr\vert\bigr] = N
\p(Y_1 \ge\sqrt{2}\gamma\log N) \le C N^{\mathcal{E}^{(\vec
{\sigma}, \alpha)}(\gamma)},
\]
where $\mathcal{E}^{(\vec{\sigma},\alpha)}(\gamma)=1- \gamma^2/V_{12}$.
Therefore, by Markov's inequality, for any $\rho>0$,
\[
\p\bigl( \bigl\vert\mathcal{H}_N^{ Y}(\gamma)\bigr \vert\ge
N^{\mathcal{E}^{(\vec
{\sigma}, \alpha)}(\gamma) + \rho} \bigr) \le C N^{-\rho}\to0,\qquad N
\to0.
\]
In case 2, if $0<\gamma< V_{12}/\sigma_1=: \gamma_{\mathrm{crit}}$ the same
argument gives the correct bound.

It remains to bound the case $\gamma\ge\gamma_{\mathrm{crit}}$.
The argument is essentially an explicit comparison with a 2-level GREM.
For the scale $\alpha$, define
\[
\mathcal{H}_{N^\alpha}^{ Y}(\gamma):=\bigl\{x \in\mathcal
X_{\varepsilon
^{\alpha}}\dvtx Y_x(\alpha) \geq\sqrt{2} \gamma\log N\bigr
\}, \qquad\mathcal{E}_1(\gamma):= \alpha- \frac{\gamma^2}{\sigma_1^2
\alpha}.
\]
A first moment computation yields, for any $0< \gamma_1<\sigma_1
\alpha
$ and any $\rho>0$,
%
%
\begin{equation}
\label{eq:controlscalealpha} \p\bigl( \bigl\vert\mathcal{H}_{N^\alpha
}^{ Y}(
\gamma_1) \bigr\vert\ge N^{\mathcal
{E}_1(\gamma_1) + \rho} \bigr) \le C N^{-\rho} \to0
,\qquad  N \to0.
\end{equation}
Similarly, a union bound gives
%
%
\begin{equation}
\label{eq:boundmax} \p\Bigl(\max_{x \in\mathcal X_{\varepsilon
^{\alpha}}} Y_x(\alpha)
\geq\sqrt{2} \sigma_1 \alpha\log N\Bigr) \to0.
\end{equation}
Recall that, for any $x \in\mathcal X_{\varepsilon}$, we denote by
$\pi
(x)$ the closest point in $\mathcal X_{\varepsilon^{\alpha}}$,
hence $\Vert x -\pi(x) \Vert\le\varepsilon^{\alpha}/2$. We define for
all $N$ and $\nu>0$,
\[
A_{N,\nu}:= \bigcup_{x \in\mathcal X_{\varepsilon}} \bigl\{ \bigl\vert
Y_x(\alpha) - Y_{\pi(x)}(\alpha) \bigr\vert\ge\nu\log N \bigr\}.
\]
The parameter $\nu$ will be fixed later and will depend on $\rho$.
Using a union bound together with Lemma \ref{lem: lem12}, we obtain,
for all $\nu>0$,
%
%
\begin{equation}
\label{eq:ANnu} \p( A_{N,\nu}) \le C N \mathrm{e}^{- c (\log N)^2} \to0,\qquad N
\to0.
\end{equation}
We also consider the events giving the log-number of high points at
scale $\alpha$.
Precisely, we divide $ [0, \sigma_1 \alpha]$ in intervals of
size $\sigma_1 \alpha/M$ where $M$ will be fixed later.
Define $\eta_i:=i \sigma_1 \alpha/M$, for $0 \le i \le M$ and
\[
I^{(i)}:= [\sqrt{2} \eta_{i-1} \log N;\sqrt{2}
\eta_{i} \log N ], \qquad 1 \le i \le M.
\]
By \eqref{eq:controlscalealpha}, the events
\[
B_{N,i}:= \bigl\{\bigl\vert\mathcal{H}_{N^\alpha}^{ Y}(
\eta_{i-1}) \bigr\vert\ge N^{\mathcal{E}_1(\eta_{i-1}) + \rho/2} \bigr
\}, \qquad 1 \le i \le M
\]
are such that
%
%
\begin{equation}
\label{eq:BNi} \p\Biggl(\bigcup_{i=1}^{M}
B_{N,i} \Biggr) \to0,\qquad  N \to0.
\end{equation}
Therefore, by \eqref{eq:ANnu} and \eqref{eq:BNi}, we are reduced to estimate
\[
\p\Biggl( \bigl\{\bigl\vert\mathcal{H}_N^{ Y}(\gamma) \bigr\vert
\ge N^{\mathcal
{E}^{(\vec{\sigma}, \alpha)}(\gamma) + \rho}\bigr\} \cap A_{N,\nu
}^{ c}\cap\bigcap
_{i=1}^M B_{N,i}^{ c}
\Biggr),
\]
which is smaller than
%
%
\begin{equation}
\label{eqn: to bound} \frac{1}{N^{\mathcal{E}^{(\vec{\sigma},
\alpha)}(\gamma) + \rho}}\E\Biggl[\bigl\vert\mathcal{H}_N^{ Y}(
\gamma) \bigr\vert; A_{N,\nu}^{ c}, \bigcap
_{i=1}^M B_{N,i}^{ c} \Biggr].
\end{equation}
We split the set $\mathcal{H}_{ N}^{ Y}(\gamma)$ into the possible
value of the field at scale $\alpha$
\begin{eqnarray*}
\mathcal{H}^{(i)}_N(\gamma)&:=& \bigl\{ x \in\mathcal
X_{\varepsilon} \dvtx Y_x \ge\sqrt{2} \gamma\log N;
Y_{\pi(x)}(\alpha) \in I^{(i)} \bigr\},\qquad  1 \le i \le M,
\\
\mathcal{H}^{(0)}_N(\gamma)&:=& \bigl\{ x \in\mathcal
X_{\varepsilon} \dvtx Y_x \ge\sqrt{2} \gamma\log N;
Y_{\pi(x)}(\alpha) \le0 \bigr\}.
\end{eqnarray*}
The term in \eqref{eqn: to bound} can then be bounded above by
\[
\frac{1}{N^{\mathcal{E}^{(\vec{\sigma}, \alpha)}(\gamma) + \rho
}} \sum_{i=0}^M \E
\bigl[\bigl\vert\mathcal{H}^{(i)}_N(\gamma) \bigr\vert;
A_{N,\nu}^{
c}\cap B_{N,i}^{ c} \bigr].
\]
If $ 0 \le\gamma\le\gamma_{\mathrm{max}}$, note that $\mathcal{E}^{(\vec
{\sigma},\alpha)}(\gamma)$ satisfies $\mathcal{E}^{(\vec{\sigma
},\alpha
)}(\gamma)= \max_{0 \le\eta\le\sigma_1 \alpha} Q(\eta)$ where
\[
Q(\eta):=1-\frac{\eta^2}{\sigma_1^2 \alpha} -\frac{(\gamma-\eta
)^2}{\sigma_2^2 (1-\alpha)}.
\]
Moreover, if $ \gamma_{\mathrm{crit}} \le\gamma\le\gamma_{\mathrm{max}}$, the maximum
is attained at $\eta=\sigma_1\alpha$, thus $Q(\eta) \le\mathcal
{E}^{(\vec{\sigma},\alpha)}(\gamma)$ for all $\eta\in[0,\sigma_1
\alpha]$. For $1\leq i\leq M$, one gets
\begin{eqnarray*}
&&\E\bigl[\bigl\vert\mathcal{H}^{(i)}_N(\gamma) \bigr\vert;
A_{N,\nu}^{ c} \cap B_{N,i}^{ c} \bigr] \\
&&\qquad=
\E\biggl[ \sum_{x\in\mathcal X_\varepsilon} \mathbf{1}_{\{Y_x\geq
\sqrt
{2}\gamma\log N, Y_{\pi(x)}(\alpha)\in I^{(i)}\}};
A_{N,\nu}^{ c} \cap B_{N,i}^{ c} \biggr]
\\
&&\qquad\leq\E\biggl[ \sum_{x\in\mathcal X_\varepsilon} \mathbf{1}_{\{
Y_x-Y_x(\alpha)\geq\sqrt{2}(\gamma-\eta_i-\nu) \log N, Y_{\pi
(x)}(\alpha)\geq\sqrt{2}\eta_{i-1}\log N \}};
B_{N,i}^{ c} \biggr]
\\
&&\qquad \leq CN^{\mathcal{E}_1(\eta_{i-1}) + \rho/2} N^{1-\alpha
}N^{-
{(\gamma-\eta_i-\nu)^2}/{(\sigma_2^2(1-\alpha))}}
\\
&&\qquad= C N^{\rho/2}N^{1- {(\eta_{i-1})^2}/{(\sigma_1^2\alpha)}-
{(\gamma-\eta_i-\nu)^2}/{(\sigma_2^2(1-\alpha))}},
\end{eqnarray*}
where the last inequality follows by the definition of $B_{N,i}$ the
independence of the field at different scales and a Gaussian estimate.
Since $Q(\eta) \le\mathcal{E}^{(\vec{\sigma},\alpha)}(\gamma)$
for all
$\eta\in[0,\sigma_1 \alpha]$, the last term is smaller
than $CN^{\mathcal{E}^{(\vec{\sigma}, \alpha)}(\gamma) + 3\rho
/4}$ by
taking $\nu$ small enough and $M$ large enough, but fixed.
For $i=0$, a similar argument gives also the bound $CN^{\mathcal
{E}^{(\vec{\sigma}, \alpha)}(\gamma) + \rho/2}$.
Putting this back in \eqref{eqn: to bound} shows that the term goes to~$0$ as $N\to\infty$ as desired.

\subsubsection{Proof of the lower bound in Proposition \texorpdfstring{\protect\ref{prop:perturbed-highpoints}}{3.2}}
The proof of the lower bound is two-step recursion.
Two lemmas are needed. The first is a generalization of the lower bound
in Daviaud's theorem; see Theorem \ref{thm:main2} or \cite{daviaud}.

%
\begin{lemma}
\label{lem:recurrence}
Let $0 <\alpha'<\alpha'' \le1$. Suppose that the parameter $\sigma$ is
constant on the strip $[0,1]_\sim\times[\varepsilon^{\alpha''},
\varepsilon^{\alpha'}]$,
and that the event
\[
\Xi:= \bigl\{\#\bigl\{x\in\mathcal{X}_{\varepsilon^{ \alpha
'}}\dvtx Y_x
\bigl( \alpha'\bigr) \geq\sqrt{2} \gamma' \log N
\bigr\} \ge N^{\mathcal{E}'} \bigr\}
\]
is such that
\[
\p\bigl(\Xi^c\bigr) \le\exp\bigl\{-c' (\log
N)^2 \bigr\}
\]
for some $\gamma' \ge0$, $\mathcal{E}'>0$ and $c'>0$.

Let
\[
\mathcal{E}(\gamma):= \mathcal{E}' + \bigl(\alpha''-
\alpha'\bigr) - \frac{(\gamma
-\gamma')^2}{ \sigma^2 (\alpha''-\alpha')}>0.
\]
Then, for any $\gamma''$ such that $\mathcal{E}(\gamma'')>0$ and any
$\mathcal{E}< \mathcal{E}(\gamma'')$, there exists $c$
such that
\[
\p\bigl(\#\bigl\{x\in\mathcal{X}_{\varepsilon^{ \alpha''}}\dvtx
Y_x \bigl(
\alpha''\bigr) \geq\sqrt{2} \gamma''
\log N \bigr\} \le N^{\mathcal{E}} \bigr) \le\exp\bigl\{ -c (\log
N)^2 \bigr\}.
\]
\end{lemma}
We stress that $\gamma''$ may be such that $\mathcal{E}(\gamma
'')<\mathcal{E}'$.
The second lemma, which follows, serves as the starting point of the
recursion and is analogous to Lemma 8 in \cite{bolthausen-deuschel-giacomin}.
%
%
\begin{lemma}
\label{lem:init}
For any $\alpha_0$ such that $0<\alpha_0<\alpha$, there exists
$\mathcal E_0=\mathcal E_0(\alpha_0)>0$ and $c=c(\alpha_0)$ such that
\[
\p\bigl(\#\bigl\{x \in\mathcal{X}_{\varepsilon^{\alpha_0}}\dvtx
Y_x(\alpha
_0)\geq0\bigr\}\leq N^{\mathcal E_0} \bigr)\leq\exp\bigl\{-c
(\log
N)^2 \bigr\}.
\]
\end{lemma}
We first conclude the proof of the lower bound in Proposition \ref
{prop:perturbed-highpoints} using the two above lemmas.

\begin{pf*}{Proof of the lower bound of Proposition \ref
{prop:perturbed-highpoints}}
Let $\gamma$ such that $0< \gamma< \gamma_{\mathrm{max}}$.
Choose $\mathcal{E}$ such that $\mathcal{E}<\mathcal{E}^{(\vec
{\sigma},
\alpha)}(\gamma)$. It will be shown that for some $c>0$
%
%
\begin{equation}
\label{eqn: lower-} \p\bigl(\bigl\vert\mathcal{H}_N^{Y}(\gamma)
\bigr\vert\leq N^{\varepsilon} \bigr)\leq\exp\bigl\{-c (\log N)^2
\bigr
\}.
\end{equation}

By Lemma \ref{lem:init}, for $\alpha_0<\alpha$ arbitrarily close to
$0$, there exists $\mathcal E_0=\mathcal E_0(\alpha_0)>0$ and
$c_0=c_0(\alpha_0)>0$, such that
%
%
\begin{equation}
\label{eqn: lower} \p\bigl(\#\bigl\{x\in\mathcal{X}_{\varepsilon
^{ \alpha_0}}\dvtx
Y_x ( \alpha_0) \geq0 \bigr\} \leq N^{\mathcal{E}_0}
\bigr)\leq\exp\bigl\{-c_0 (\log N)^2 \bigr\}.
\end{equation}
Observe that we have $0 \le\mathcal E_0 \le\alpha_0$. Moreover, let
%
%
\begin{equation}
\label{eqn: E_1} \mathcal{E}_1(\gamma_1):=
\mathcal{E}_0 + (\alpha-\alpha_0) - \frac
{\gamma_1^2}{ \sigma_1^2 (\alpha-\alpha_0)}
.
\end{equation}
Lemma \ref{lem:recurrence} is applied from $\alpha_0$ to $\alpha$. For
any $\gamma_1$ with $\mathcal{E}_1(\gamma_1)>0$ and any
$\mathcal{E}_1<\mathcal{E}_1(\gamma_1)$, there exists $c_1>0$ such that
\[
\p\bigl(\#\bigl\{x\in\mathcal{X}_{\varepsilon^{ \alpha}}\dvtx
Y_x (
\alpha) \geq\sqrt{2}\gamma_1\log N \bigr\} \leq N^{\mathcal{E}_1}
\bigr)\leq\exp\bigl\{ -c_1 (\log N)^2 \bigr\}.
\]
Therefore, Lemma \ref{lem:recurrence} can be applied from $\alpha$ to
$1$ for any $\gamma_1$ with $\mathcal{E}_1(\gamma_1)>0$. Define similarly
%
%
\begin{equation}
\label{eqn: E_2} \mathcal{E}_2(\gamma_1,
\gamma_2):= \mathcal{E}_1(\gamma_1) + (1-
\alpha) - \frac{(\gamma_2-\gamma_1)^2}{ \sigma_2^2 (1-\alpha)}.
\end{equation}
Then, for any $\gamma_2$ with $\mathcal{E}_2(\gamma_1,\gamma_2)>0$
and $\mathcal{E}_2<\mathcal{E}_2(\gamma_1,\gamma_2)$, there exists
$c_2>0$ such that
%
%
\begin{equation}
\label{eqn: lower prob} \p\bigl(\#\{x\in\mathcal{X}_{\varepsilon
}\dvtx
Y_x \geq\sqrt{2}\gamma_2\log N \} \leq
N^{\mathcal{E}_2} \bigr)\leq\exp\bigl\{-c_2 (\log N)^2
\bigr\}.
\end{equation}
Recalling that $0 \le\mathcal E_0 \le\alpha_0$, equation \eqref{eqn:
lower-} follows from \eqref{eqn: lower prob} if it is proved that
$\lim_{\alpha_0 \to0}\mathcal{E}_2(\gamma_1,\gamma) = \mathcal
{E}^{(\vec
{\sigma},{\alpha})}(\gamma)$ for an appropriate choice of $\gamma_1$
[in particular such that $\mathcal{E}_1(\gamma_1)>0$].
It is easily verified that, for a given $\gamma$, the quantity
$\mathcal
{E}_2(\gamma_1,\gamma)$ is maximized at
\[
\gamma_1^*= \gamma\frac{\sigma_1^2(\alpha-\alpha
_0)}{V_{12}-\sigma
_1^2\alpha_0}.
\]
Plugging these back in \eqref{eqn: E_1} shows that $\mathcal
{E}_1(\gamma
_1^*)>0$ provided that
\[
\gamma< \frac{V_{12}}{\sigma_1}=:\gamma_{\mathrm{crit}},
\]
with $\alpha_0$ small enough (depending on $\gamma$).
Furthermore, since
\[
\mathcal{E}_2\bigl(\gamma_1^*,\gamma\bigr)=
\mathcal{E}_0 +(1-\alpha_0)- \frac
{\gamma^2}{V_{12}-\sigma_1^2\alpha_0},
\]
we obtain
$
\lim_{\alpha_0 \to0} \mathcal{E}_2(\gamma_1^*,\gamma) = \mathcal
{E}^{(\vec{\sigma},{\alpha})}(\gamma)$,
which completes the proof in the case $0<\gamma< \gamma_{\mathrm{crit}}$.

If $\gamma_{\mathrm{crit}} \leq\gamma< \gamma_{\mathrm{max}}$, the condition
$\mathcal
{E}_1(\gamma_1^*)>0$ will be violated as $\alpha_0$ goes to zero.
In this case, for $\nu>0$, pick $\gamma_1^{**}= \sigma_1 \alpha-\nu$
such that $\mathcal{E}_1(\gamma_1^{**})> 0$. The first term in
$\gamma
_1^{**}$ corresponds to $\gamma_1^*$ evaluated at $\gamma_{\mathrm{crit}}$ for
$\alpha_0=0$.
In particular, $\lim_{\alpha_0\to0, \nu\to0} \mathcal{E}_1(\gamma
_1^{**})= 0$. From \eqref{eqn: E_2}, this shows that
\[
\lim_{\alpha_0\to0, \nu\to0} \mathcal{E}_2\bigl(
\gamma_1^{**},\gamma\bigr)= (1-\alpha) - \frac{(\gamma-\sigma
_1\alpha)^2}{\sigma_2^2(1-\alpha)}=
\mathcal{E}^{(\vec{\sigma},{\alpha})}(\gamma).
\]
Note that $\mathcal{E}^{(\vec{\sigma},{\alpha})}(\gamma)$ is strictly
positive if and only if $\gamma< \sigma_1\alpha+\sigma_2(1-\alpha
)=\gamma_{\mathrm{max}}$.
This concludes the proof of \eqref{eqn: lower-}.
\end{pf*}

\begin{pf*}{Proof of Lemma \ref{lem:recurrence}}
Let $\gamma''$ such that $\mathcal{E}(\gamma'')>0$ and $\mathcal{E}$
such that $0<\mathcal{E}<\mathcal{E}(\gamma'')$.
Pick $\overline{\gamma}>\gamma''$ such that
%
%
\begin{equation}
\label{eq:rec1} \mathcal{E}(\overline{\gamma})> \mathcal{E}>0.
\end{equation}

Since $\overline\gamma>\gamma''$, there exists $\varsigma\in(0,1)$
such that
%
%
\begin{equation}
\label{eq:rec2} \overline{\gamma} (1-\varsigma) \ge\gamma''.
\end{equation}
For $K \in{\mathbb N}$ (which will be fixed later), we set
\begin{eqnarray*}
\eta_\ell&:=& \alpha' +\frac{\ell-1}{K} \bigl(
\alpha''-\alpha'\bigr),\qquad 1 \le\ell\le
K+1,
\\
\lambda_\ell&:=&\gamma'+ \frac{\ell-1}{K} \bigl(
\overline{\gamma}-\gamma'\bigr) (1-\varsigma),\qquad  1 \le\ell\le K+1.
\end{eqnarray*}
Observe that the $\eta_\ell$'s and the $\lambda_\ell$'s satisfy
$
\eta_1= \alpha'<\eta_2<\cdots<\eta_{K}<\eta_{K+1}= \alpha''$,
and
$
\lambda_1=\gamma'<\lambda_2<\cdots<\lambda_{K}<\lambda_{K+1}=
(1-\varsigma)\overline{\gamma} + \varsigma\gamma'$.
Consider the sets $\mathcal{A}_\ell$ given by
\[
\mathcal{A}_\ell:= \bigl\{ \underline{x}^{(\ell)}=(x_1,
\ldots,x_\ell) \dvtx x_i \in\mathcal{X}_{2 \varepsilon^{\eta_i}},
\forall1 \le i \le\ell\mbox{ and }  \Vert x_{i+1}-x_i \Vert
\le\varepsilon^{\eta_i}/2 \bigr\}
\]
for $1 \le\ell\le K+1$. Note that only half of the $x_i$'s in
$\mathcal{X}_{\varepsilon^{\eta_i}}$'s are considered.
Also, to each $x_i$ we consider the points $x_{i+1}$ in $\mathcal
{X}_{2\varepsilon^{\eta_{i+1}}}$
that are close to $x_i$. By analogy with a branching process, these
points can be thought of as the \textit{children} of $x_{i}$.
The reason for these two choices is that the cones corresponding to the
variables $Y_{x_{i+1}}(\eta_{i+1})$ and $Y_{x'_{i+1}}(\eta_{i+1})$ do
not intersect below the line $y=\varepsilon^{\eta_i}$
if $x_i\neq x_i'$; see Figure~\ref{fig:intersect-daviaud}.

%
\begin{figure}

\includegraphics{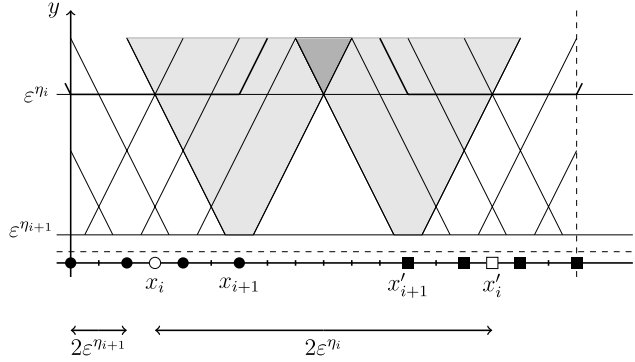}

\caption{Approximation by a tree-like structure. The black circles
symbolize the children of the white circle, while the black squares
symbolize the children of the white square.}
\label{fig:intersect-daviaud}
\end{figure}

Now consider, the sets of high points of $\mathcal{A}_\ell$,
\[
A_\ell:= \bigl\{ \underline{x}^{(\ell)} \in
\mathcal{A}_\ell\dvtx Y_{x_i} (\eta_i) \ge
\sqrt{2} \lambda_i \log N, \forall1 \le i \le\ell\bigr\}, \qquad 1 \le
\ell\le K+1
\]
and
\[
B_\ell:= \{ \# A_\ell\ge n_\ell\}, \qquad 1 \le\ell
\le K+1,
\]
where
%
%
\begin{equation}
\label{eqn: nl} \qquad n_\ell:= N^{\mathcal{E}'+{(\ell-1)}/{K} ((\alpha
''-\alpha')-
({(\overline{\gamma}-\gamma')^2}/{(\sigma^2 (\alpha''-\alpha'))}) )},
\qquad 1 \le\ell\le K+1,
\end{equation}
such that
$
N^{\mathcal{E}'}=n_1$ and $n_{K+1}=N^{\mathcal{E}(\overline\gamma)}$.
Furthermore, with these definitions and the choice of $\overline\gamma
$ in \eqref{eq:rec2} and \eqref{eq:rec1}, we have for large $N$
\begin{eqnarray*}
B_{K+1}&=& \{ \# A_{K+1} > n_{{K+1}} \}
\\
&&{}\subset\bigl\{ \#\bigl\{x\in\mathcal{X}_{\varepsilon^{\alpha
''}}\dvtx Y_x
\bigl(\alpha''\bigr) \geq\sqrt{2}\bigl((1-\varsigma)
\overline\gamma+ \varsigma\gamma'\bigr)\log N\bigr\}>
N^{\mathcal{E}(\overline\gamma)} \bigr\}
\\
&&{}\subset\bigl\{ \#\bigl\{x\in\mathcal{X}_{\varepsilon^{\alpha
''}}\dvtx Y_x
\bigl(\alpha''\bigr) \geq\sqrt{2}\gamma''
\log N\bigr\}> N^{\mathcal E} \bigr\}.
\end{eqnarray*}
It is thus sufficient to find a bound for $\p(B_{K+1}^c)$ to prove the
lemma. For events $C_\ell$ to be defined in \eqref{eqn: C}, we use the
elementary bound
$\p(B_{K+1}^c) \le\p(B_{K+1}^c \cap B_{K} \cap C_{K}^c) + \p(C_{K}) +
\p(B_{K}^c)$ which applied recursively gives
%
%
\begin{equation}
\label{eqn: sum events} \p\bigl(B_{K+1}^c\bigr) \le\sum
_{\ell=2}^{K+1} \bigl(\p\bigl(B_{\ell}^c
\cap B_{\ell
-1} \cap C_{\ell-1}^c\bigr) +
\p(C_{\ell-1}) \bigr) + \p\bigl(B_{1}^c\bigr).
\end{equation}
The last term has the correct bound by assumption. It remains to bound
the ones appearing in the sum.

On the event $B_\ell$, there exist at least $n_\ell$ high $\ell
$-branches $\underline{x}^{(\ell)}=(x_1,\ldots,x_\ell)$, these are
branches that satisfy $Y_{x_i} (\eta_i) \ge\sqrt{2} \lambda_i \log
N $
for $1 \le i \le\ell$.
Select the first $n_\ell$ such $\ell$-branches, and denote them by
$
\underline{x}^{(\ell)}_j=(x_{j,1},\ldots,x_{j,\ell})$,
for all $1 \le j \le n_\ell$.
Consider the set $\mathcal{A}_{j,\ell}$, the children of $x_{j,\ell}$
at level $\eta_{\ell+1}$:
$
\mathcal{A}_{j,\ell}:= \{ x \in\mathcal{X}_{2 \varepsilon^{\eta
_{\ell+1}}} \dvtx\Vert x- x_{j,\ell}\Vert\le\varepsilon^{\eta
_\ell}/2
\}
$.
It holds
\begin{eqnarray*}
B_{\ell} \cap B_{\ell+1}^c &\subset&
B_{\ell} \cap\Biggl\{\sum_{j=1}^{n_\ell}
\sum_{x \in\mathcal{A}_{j,\ell}} \mathbf{1}_{ \{
Y_{x} (\eta_{\ell+1}) - Y_{x_{j,\ell}} (\eta_{\ell}) \ge\sqrt{2}
({(\overline{\gamma}-\gamma') (1-\varsigma)}/{K}) \log N \}} \le
n_{\ell+1} \Biggr\}
\\
&\subset& B_{\ell} \cap\Biggl\{\sum_{j=1}^{n_\ell}
\zeta_j \le\frac{2
n_{\ell+1}}{N^{(\alpha''-\alpha')/K}} \Biggr\},
\end{eqnarray*}
where
%
%
\begin{equation}
\label{eqn: zeta} \zeta_j:= \frac{1}{\vert\mathcal{A}_{j,\ell}
\vert} \sum
_{x \in
\mathcal{A}_{j,\ell}} \mathbf{1}_{ \{ Y_{x} (\eta_{\ell+1}) -
Y_{x_{j,\ell}} (\eta_{\ell}) \ge\sqrt{2} ({(\overline{\gamma
}-\gamma') (1-\varsigma)}/{K}) \log N \}},
\end{equation}
and $\vert\mathcal{A}_{j,\ell} \vert= N^{(\alpha''-\alpha')/K}/2$.
A crucial point is that $Y_{x_{j,\ell}}(\eta_\ell)$ is not equal to
$Y_x(\eta_\ell)$ since $x\neq x_{j,\ell}$ in general.
However, it turns out that their value must
be very close since the variance of the difference is essentially a
constant due to the logarithmic correlations. Precisely,
let
%
%
\begin{eqnarray}
\label{eqn: C} C_{\ell}&:=& \bigcup_{\underline{x}^{(\ell)} \in
\mathcal{A}_\ell}
\mathop{\bigcup_{x \in\mathcal{X}_{2 \varepsilon^{\eta_{\ell
+1}}} \dvtx}}_{
\Vert x- x_{\ell}\Vert\le\varepsilon^{\eta_\ell}/2} \biggl\{ \bigl|
Y_{\underline{x}^\ell} (\eta_{\ell}) - Y_{x} (
\eta_{\ell
}) \bigr|
\nonumber
\\[-8pt]
\\[-8pt]
\nonumber
&&\hspace*{90pt}{}\ge\sqrt{2} \nu\frac{(\overline{\gamma}-\gamma')
(1-\varsigma)}{K} \log N \biggr\}
\end{eqnarray}
for $\nu>0$ which is fixed and will be chosen small later.
By Lemma \ref{lem: lem12} of the \hyperref[app]{Appendix},
$
\operatorname{Var}( Y_{x} (\eta_{\ell}) - Y_{x'} (\eta_{\ell}) ) \le\max\{
\sigma^2_1,\sigma_2^2\}<\infty$,
for every $1 \le\ell\le K$, and any $x \in\mathcal{X}_{2
\varepsilon
^{\eta_{\ell}}}$, $x' \in\mathcal{X}_{2 \varepsilon^{\eta_{\ell+1}}}$
such that $ \Vert x'-x\Vert\le\varepsilon^{\eta_\ell}/2$.
Therefore, a Gaussian estimate (see Lemma \ref{lem: gaussian}),
together with the union-bound give
%
%
\begin{equation}
\label{eq:rec4} \p(C_{\ell}) \le\exp\bigl\{-d (\log N)^2
\bigr\}
\end{equation}
for all $1 \le\ell\le K$ and some $d>0$.

It remains to bound the first term appearing in the sum of \eqref{eqn:
sum events}.
On $C_{\ell}^c$, $Y_{x_{j,\ell}} (\eta_{\ell})$ can be replaced by $Y_x
(\eta_{\ell})$ in \eqref{eqn: zeta}, making a small error that depends
on $\nu$. Namely, one has $\zeta_j \ge\tilde{\zeta}_j$, where
\[
\tilde{\zeta}_j:= \frac{1}{\vert\mathcal{A}_{j,\ell} \vert} \sum
_{x
\in\mathcal{A}_{j,\ell}} \mathbf{1}_{ \{ Y_{x} (\eta_{\ell+1}) -
Y_{x} (\eta_{\ell}) \ge\sqrt{2} (1+\nu)( {(\overline{\gamma
}-\gamma
')(1-\varsigma)}/{K}) \log N \}}.
\]
Note that conditionally on $\mathcal{F} _{\varepsilon^{\eta_\ell}}$,
the $\tilde{\zeta}_j$'s are i.i.d. Moreover, since the $ \tilde
{\zeta
}_j$'s are independent of $\mathcal{F} _{\varepsilon^{\eta_\ell}}$,
they are also independent of each other. Lemma \ref{lem: LD} of the
\hyperref[app]{Appendix} guarantees that the sum of the $\tilde\zeta
_j$ cannot be too low.
Observe that
\[
{\mathbb E}[ \tilde{\zeta}_j ] = \p\biggl(z \ge\sqrt{2} (1+\nu)
\frac
{(\overline{\gamma}-\gamma')(1-\varsigma)}{K} \log N \biggr),
\]
where $z$ is a centered Gaussian with variance
$
\sigma^2 \log( \frac{\varepsilon^{\eta_{\ell}}}{\varepsilon
^{\eta
_{\ell+1}}} ) = \sigma^2 \frac{(\alpha''-\alpha')}{K} \log N$.
By a Gaussian estimate, Lemma \ref{lem: gaussian},
\[
{\mathbb E}[ \tilde{\zeta}_j ] \ge\exp\biggl\{ -\frac{1}{K}
\frac
{(1+2\nu)^2 (\overline{\gamma}-\gamma')^2 (1-\varsigma)^2}{\sigma^2
(\alpha''-\alpha')} \log N \biggr\},
\]
where $(1+\nu)$ has been replaced by $(1+2\nu)$ to absorb the
$1/\sqrt{\log N}$ term in front of the exponential.
Consequently, using elementary manipulations,
\begin{eqnarray*}
&&B_{\ell+1}^c\cap B_\ell\cap C_\ell^c\\
&&\qquad\subset \Biggl\{ \sum_{j=1}^{ n_\ell} \bigl(
\tilde{\zeta}_j - {\mathbb E}[\tilde{\zeta}_j ] \bigr)\\
&&\hspace*{38pt}{} \le
\frac{2n_{\ell+1}}{N^{(\alpha''-\alpha')/K}} - n_\ell N^{-
({1}/{K})({(1+2\nu)^2 (\overline{\gamma}-\gamma')^2 (1-\varsigma
)^2}/{(\sigma^2 (\alpha''-\alpha'))})} \Biggr\}
\\
&&\qquad\subset \Biggl\{ \Biggl\llvert\sum_{j=1}^{n_\ell}
\bigl(\tilde{\zeta}_j - {\mathbb E}[\tilde{\zeta}_j ] \bigr)
\Biggr\rrvert\ge\frac{1}{2} n_\ell N^{-({1}/{K})({(1+2\nu)^2
(\overline{\gamma}-\gamma')^2 (1-\varsigma)^2}/{(\sigma^2 (\alpha
''-\alpha
'))})} \Biggr\},
\end{eqnarray*}
provided
\[
\frac{1}{K}\frac{(1+2\nu)^2 (\overline{\gamma}-\gamma')^2
(1-\varsigma
)^2}{\sigma^2 (\alpha''-\alpha')} < \frac{1}{K} \frac{(\overline
{\gamma
}-\gamma')^2}{\sigma^2 (\alpha''-\alpha')},
\]
that is %
%
\begin{equation}
\label{eq:rec5} (1+2\nu) (1-\varsigma) <1.
\end{equation}
Fix $\nu$ small enough such that (\ref{eq:rec5}) is satisfied.
Write for short
\[
\mu:=\frac{1}{K}\frac{(1+2\nu)^2 (\overline{\gamma}-\gamma')^2
(1-\varsigma)^2}{\sigma^2 (\alpha''-\alpha')}.
\]
Then, taking $n=n_\ell$ and $t=n_\ell N^{-\mu}$ in Lemma \ref{lem: LD},
we get
\begin{eqnarray}
\nonumber
\p\bigl(B_{\ell+1}^c\cap B_\ell\cap
C_\ell^c \bigr) &\le& 2 \exp\biggl\{ \frac
{n_\ell^2 N^{-2 \mu}}{2n_\ell+ ({2}/{3}) n_\ell N^{-\mu}}
\biggr\}. \label{eq:rec7}
\end{eqnarray}
By the form of $n_{\ell}$ in \eqref{eqn: nl}, $K$ can be taken large
enough so that $n_{\ell} N^{-2\mu}>N^\delta$ for some $\delta>0$ and
all $\ell=1,\ldots,K+1$.
This completes the proof of the lemma.
\end{pf*}

\begin{pf*}{Proof of Lemma \ref{lem:init}}
Take $\alpha'<\alpha_0$ in such a way that $\mathcal{X}_{\varepsilon
^{\alpha'}} \subset\mathcal{X}_{\varepsilon^{\alpha_0}}$.
Consider the set
\[
\Lambda:= \bigl\{x\in\mathcal{X}_{\varepsilon^{\alpha'}}\dvtx
Y_x\bigl(
\alpha'\bigr)\geq-\sigma_1\bigl(\alpha_0-
\alpha'\bigr)\log N \bigr\},
\]
and the event
\[
A=A_\delta:= \bigl\{|\Lambda|\geq N^\delta\bigr\},\qquad \delta>0.
\]
The parameters ${\mathcal E_0}$, $\delta$ and $\alpha'$ will be chosen
later as a function of $\alpha_0$.
By splitting the probability on the event $A$,
\begin{eqnarray*}
&&\p\bigl(\#\bigl\{x\in\mathcal{X}_{\varepsilon^{\alpha_0}}\dvtx
Y_x(\alpha
_0)\geq0\bigr\}\leq N^{\mathcal E_0} \bigr)
\\
&&\qquad\leq\p\bigl(\#\bigl\{x\in\mathcal{X}_{\varepsilon^{\alpha
_0}}\dvtx Y_x(
\alpha_0) \geq0\bigr\}\leq N^{\mathcal E_0};A \bigr)+\p
\bigl(A^c\bigr)
\\
&&\qquad\leq{\mathbb E}\bigl[ \p\bigl(\#\bigl\{x \in\Lambda\dvtx Y_x(
\alpha_0)-Y_x\bigl(\alpha'\bigr) \geq
\sigma_1 \bigl(\alpha_0-\alpha'\bigr)\log
N \bigr\}\leq N^{\mathcal E_0} | {\mathcal{F}}_{\varepsilon^{\alpha
'}} \bigr);A \bigr]\\
&&\qquad\quad{}+\p
\bigl(A^c\bigr),
\end{eqnarray*}
where the second inequality is obtained by restricting to the set
$\Lambda\subset\mathcal{X}_{\varepsilon^{\alpha_0}}$.

First we prove that the definition of $A$ yields a super-exponential
decay of the first term for ${\mathcal E_0}$ and $\delta$ depending on
$\alpha_0-\alpha'$.
The variables $Y_x(\alpha_0)-Y_x(\alpha')$, $x \in\mathcal
{X}_{\varepsilon^{\alpha'}}$, are i.i.d. Gaussians of variance
$\sigma
_1^2(\alpha_0-\alpha')\log N$.
Write for simplicity $(z_i, i=1,\ldots,N^\delta)$ for i.i.d. Gaussians
random variables with variance $\sigma_1^2(\alpha_0-\alpha')\log N$.
A Gaussian estimate (see Lemma \ref{lem: gaussian}) implies
\[
\p\bigl(z_i\geq\sigma_1\bigl(\alpha_0-
\alpha'\bigr)\log N \bigr)\geq\frac
{1}{2} \frac{\ee^{-({1}/{2})(\alpha_0-\alpha')\log N}}{\sqrt
{(\alpha
_0-\alpha')\log N}}
\geq\ee^{-({2}/{3})(\alpha_0-\alpha')\log N}.
\]
Therefore
\begin{eqnarray*}
& &{\mathbb E}\bigl[ \p\bigl(\#\bigl\{x \in\Lambda\dvtx Y_x(
\alpha_0)-Y_x\bigl(\alpha'\bigr) \geq
\sigma_1 \bigl(\alpha_0-\alpha'\bigr)\log
N \bigr\}\leq N^{\mathcal E_0} | {\mathcal{F}}_{\varepsilon^{\alpha
'}} \bigr);A \bigr]
\\
&&\qquad\leq\p\Biggl(\sum_{i=1}^{N^\delta} \bigl(
\mathbf{1}_{\{z_i\geq\sigma
_1 (\alpha_0-\alpha')\log N\}}-\p\bigl(z_i\geq\sigma_1
\bigl(\alpha_0-\alpha'\bigr)\log N \bigr) \bigr)\\
&&\hspace*{166pt}{} \leq
N^{\mathcal E_0}-N^{\delta-({2}/{3})(\alpha_0-\alpha')} \Biggr).
\end{eqnarray*}
Lemma \ref{lem: LD} in the \hyperref[app]{Appendix} gives a
super-exponential decay of
the above probability for the choice $ \delta>\frac{4}{3}(\alpha
_0-\alpha')$ and ${\mathcal E_0}-\delta+\frac{2}{3}(\alpha_0-\alpha
')<0$, for example,
$\delta=2(\alpha_0-\alpha')$ and ${\mathcal E_0}=\alpha_0-\alpha'$.

It remains to show that $\p(A^c)$ has super-exponential decay.
We have
\[
\p\bigl(A^c\bigr)\leq P\Bigl(A^c, \max
_{x \in\mathcal{X}_{\varepsilon^{\alpha'}}} Y_x\bigl(\alpha'\bigr
)\leq(
\log N)^2\Bigr)+P\Bigl( \max_{x \in\mathcal{X}_{\varepsilon
^{\alpha'}}}
Y_x\bigl(\alpha'\bigr)> (\log N)^2
\Bigr).
\]
The second term is easily shown to have the desired decay. We focus on
the first.
On the event $A^c\cap\{ \max_{x \in\mathcal{X}_{\alpha'}}
Y_x(\alpha
')\leq(\log N)^2\}$,
%
%
\begin{eqnarray}
\label{eqn: init} %
&&\frac{1}{| \mathcal{X}_{\varepsilon^{\alpha'}}|}\sum_{x\in
\mathcal
{X}_{\varepsilon^{\alpha'}}}
\omega_{\alpha'}(x)\nonumber\\
&&\qquad=\frac{1}{| \mathcal
{X}_{\varepsilon^{\alpha'}}|}\sum_{x\in\Lambda}
\omega_{\alpha'}(x)+\frac
{1}{| \mathcal{X}_{\varepsilon^{\alpha'}}|}\sum_{x\in\Lambda^c}
\omega_{\alpha'}(x)
\\
& &\qquad\leq\frac{|\Lambda|}{| \mathcal{X}_{\varepsilon^{\alpha
'}}|}(\log N)^2 + \biggl(1-\frac{|\Lambda|}{| \mathcal
{X}_{\varepsilon^{\alpha
'}}|}
\biggr) \bigl(-\sigma_1 \bigl(\alpha_0-
\alpha'\bigr)\log N\bigr).\nonumber
\end{eqnarray}
Since $| \mathcal{X}_{\varepsilon^{\alpha'}}|=N^{\alpha'}$, it is
easily checked that for $\delta=2(\alpha_0-\alpha')<\alpha'$, the above
is smaller than $-\frac{2}{3} \sigma_1 (\alpha_0-\alpha')\log N$.
Therefore we choose $\alpha'$ such that $\alpha_0<3\alpha'/2$.
Finally the left-hand side of \eqref{eqn: init} is a Gaussian random
variable, whose variance is of order $1$. Therefore the probability
that it is smaller than $-\frac{2}{3}\sigma_1(\alpha_0-\alpha')\log N$
is super-exponentially small.
This completes the proof of the lemma.
\end{pf*}

\section{The free energy from the high points: Proof of Proposition
\texorpdfstring{\protect\ref{prop:perturbed-free-energy}}{2.1}}
\label{sect: free energy}

In this section, we compute the free energy of the perturbed models
introduced in Section~\ref{sect: perturbed}.
The free energy $f^{(\vec{\sigma}, \alpha)}_N(\beta)$ is shown to
converge in probability to the claimed expression.
The $L^1$-convergence then follows from the fact that the variables
$(f^{(\vec{\sigma}, \alpha)}_N(\beta))_{N \ge1}$ are uniformly integrable.
This is a consequence of Borell-TIS inequality.
(Another more specific approach used by Capocaccia, Cassandro and Picco
\cite{capocaccia-cassandro-picco} for the GREM models could also have
been applied here; see Section~3.1 in \cite{capocaccia-cassandro-picco}.
Indeed, we clearly have
\[
\beta\frac{\max_{x \in\mathcal X_{N}} Y_x}{\log N} \leq f^{(\vec
{\sigma}, \alpha)}_N(\beta) \leq1+\beta
\frac{\max_{x \in\mathcal
X_{N}} Y_x}{\log N}.
\]
Therefore, uniform integrability follows if it is proved that
$
\frac{1}{ ( \log N )^2 }\times\break  \E[ ( \max_{x \in\mathcal X_{N}} Y_x
)^2 ]
$
is uniformly bounded. It equals
\[
\frac{1}{ ( \log N )^2 }\E\Bigl[ \Bigl( \max_{x \in\mathcal X_{N}}
Y_x - \E\Bigl[ \max_{x \in\mathcal X_{N}} Y_x
\Bigr] \Bigr)^2 \Bigr] + \frac{1}{ (
\log N )^2 }\E\Bigl[ \max
_{x \in\mathcal X_{N}} Y_x\Bigr] ^2.
\]
The first term is bounded by the Borell-TIS inequality (see \cite
{adler-taylor}, page 50)
\[
\p\Bigl( \Bigl|\max_{x \in\mathcal X_{N}} Y_x - \E\max
_{x \in
\mathcal X_{N}} Y_x \Bigr| >r \Bigr)\leq2 \mathrm{e}^{-{r^2}/{(2V_{12} \log
N)}}\qquad
\forall r>0,
\]
which gives
\[
\E\biggl[ \biggl( \frac{\max_{x \in\mathcal X_{N}} Y_x - \E[
\max_{x
\in\mathcal X_{N}} Y_x]}{\log N} \biggr)^2 \biggr] \leq4 \int
_0^\infty r \mathrm{e}^{-{r^2}/{(2V_{12})}\log N } \,\d r.
\]
The right-hand side goes to zero for $N\to\infty$.
The term $\frac{1}{\log N}\E[ \max_{x \in\mathcal X_{N}} Y_x]$ can be
bounded uniformly by comparing with i.i.d. centered Gaussian random
variables of variance $V_{12} \log N$ and using Slepian's inequality;
see, for example, \cite{adler-taylor}, page~57.
Equivalently, one can reason as follows.
It is easily checked that the probability that the maximum be negative
decreases exponentially with $N$.
Thus to control the second term it suffices to control
\[
\frac{1}{\log N} \int_0^\infty\p\Bigl(\max
_{x \in\mathcal X_{N}} Y_x>r\Bigr) \,\d r.
\]
It suffices to split the integral in two intervals: $[0,\sqrt
{2V_{12}}\log N)$ and\break  $[\sqrt{2V_{12}}\log N, +\infty)$. The first
integral divided by $\log N$ is evidently of order $1$.
The second integral divided by $\log N$ tends to $0$ by a union bound
and a Gaussian estimate.
The almost-sure convergence is straightforward from the
$L^1$-convergence and the almost-sure self-averaging property of the
free energy
\[
\lim_{N\to\infty}\bigl\vert f^{(\vec{\sigma}, \alpha)}_N(\beta)- \E
f^{(\vec
{\sigma}, \alpha)}_N(\beta)\bigr\vert= 0 \qquad\mbox{a.s. }
\]
This is a standard consequence of concentration of measure (see \cite
{talagrand}, page~32)
since the free energy is a Lipschitz function of i.i.d. Gaussian
variables of Lipschitz constant smaller than $\beta/\sqrt{\log N}$.
(Note that the $Y_x$'s can be written as a linear combination of
i.i.d. standard Gaussians with coefficients chosen to get the correct
covariances.)

It remains to prove that the free energy $f^{(\vec{\sigma}, \alpha
)}_N(\beta)$ converges in probability to the claimed expression in
Proposition \ref{prop:perturbed-free-energy}.
For fixed $\beta>0$ and $\nu>0$, we prove that
%
%
\begin{eqnarray}
\label{eq:tm3-eq1} \lim_{N\to\infty}\p\bigl( f^{(\vec{\sigma},
\alpha)}_N(
\beta)\le f^{(\vec{\sigma}, \alpha)}(\beta) - \nu\bigr)&=&0,
\\
\lim_{N\to\infty}\p\bigl(f^{(\vec{\sigma}, \alpha)}_N(\beta)
\ge f^{(\vec{\sigma}, \alpha)}(\beta)+ \nu\bigr)&=&0. \label{eq:tm3-eq2}
\end{eqnarray}

First, we introduce some notation and give a preliminary result.
For simplicity, we will write $\mathcal{E}$ for $\mathcal{E}^{(\vec
{\sigma}, \alpha)}$ throughout the proof.
For any $M\in{\mathbb N}$, consider the partition of $ [0,\gamma_{\mathrm{max}}
]$ into $M$ intervals $ [\gamma_{i-1},\gamma_{i} [$,
where the $\gamma_i$'s are given by
\[
\gamma_i:=\frac{i}{M} \gamma_{\mathrm{max}}, \qquad i=0,1,
\ldots,M.
\]
Moreover for any $N \ge2$, any $M \in{\mathbb N}$ and any $\delta
>0$, define
the random variable
\[
K_{N,M}(i) :=\# \biggl\{x \in\mathcal{X}_N \dvtx
\frac{Y_x}{\sqrt{2} \log
N} \in[\gamma_{i-1},\gamma_{i} [ \biggr\},\qquad
1 \le i \le M,
\]
and the events
\begin{eqnarray*}
\nonumber
B_{N,M,\delta}&:=&\bigcap_{i=1}^M
\bigl\{N^{\mathcal{E}(\gamma
_{i-1})-\delta}-N^{\mathcal{E}(\gamma_i)+\delta} \le K_{N,M}(i) \le
N^{\mathcal{E}(\gamma_{i-1})+\delta} -N^{\mathcal{E}(\gamma
_{i})-\delta
} \bigr\}
\\
&&{} \cap\bigl\{ \#\{x \in\mathcal{X}_N \dvtx
Y_x \ge{\sqrt{2} \gamma_{\mathrm{max}} \log N}\} =0 \bigr\}.
\end{eqnarray*}
The next result is a straightforward consequence of Propositions~\ref
{prop:perturbed-maximum} and~\ref{prop:perturbed-highpoints}.

%
\begin{lemma}
\label{lem:thm3-Bonenv}
For any $M \in{\mathbb N}$ and any $\delta>0$, we have
\[
\lim_{N \to\infty} \p(B_{N,M,\delta} )=1.
\]
\end{lemma}

Define the continuous function
\[
P_{\beta}(\gamma):=\mathcal{E}(\gamma) + \sqrt{2}\beta\gamma\qquad
\forall
\gamma\in[0,\gamma_{\mathrm{max}} ].
\]
Using the expression of $\mathcal{E}$ in Proposition \ref
{prop:perturbed-highpoints} on the different intervals,
it is easily checked by differentiation that
%
%
\begin{equation}
\label{eq:sup=free} \max_{\gamma\in[0,\gamma_{\mathrm{max}} ]} P_{\beta
}(\gamma
)=f^{(\vec{\sigma}, \alpha)}(\beta).
\end{equation}
Furthermore, the continuity of $\gamma\mapsto P_{\beta}(\gamma)$ on
$ [0,\gamma_{\mathrm{max}} ]$ yields
\[
\max_{0 \le i \le M-1} P_{\beta}(\gamma_i)
\longrightarrow\max_{\gamma
\in[0,\gamma_{\mathrm{max}} ]} P_{\beta}(
\gamma)=f^{(\vec{\sigma},
\alpha)}(\beta), \qquad M \to\infty.
\]

Fix $M \in{\mathbb N}$ large enough and $\delta>0$ small enough, such that
%
%
\begin{eqnarray}
\label{eq:tm3-choiceM} \max_{0 \le i \le M-1} P_{\beta}(
\gamma_i) &\ge& f^{(\vec{\sigma},\vec
{\alpha})}(\beta) - \frac{\nu}{3},
\\
\label{eq:tm3-choiceM-bis} \frac{\sqrt{2}\beta}{M}&<&\frac{\nu}{3},
\\
\delta&<&\min\biggl\{ - \frac{1}{2} \max_{1 \le i \le M}\bigl\{
\mathcal{E}(\gamma_{i})-\mathcal{E}(\gamma_{i-1})\bigr\},
\frac{\nu}{3},\sqrt{2}\gamma_1\beta\biggr\}. \label{eq:tm3-choicedelta}
\end{eqnarray}
Note that for fixed $M$, $ \max_{1 \le i \le M}\{\mathcal{E}(\gamma
_{i})-\mathcal{E}(\gamma_{i-1})\}<0$ since $\gamma\mapsto\mathcal
{E}(\gamma)$ is a decreasing function on $ [0,\gamma_{\mathrm{max}} ]$.

\textit{Proof of the lower bound} (\ref{eq:tm3-eq1}).
Observe that the partition function\break  $Z^{(\vec{\sigma}, \alpha
)}_N(\beta
)$ associated with the perturbed model satisfies
$
Z^{(\vec{\sigma}, \alpha)}_N(\beta) \ge\break \sum_{i=1}^{M}K_{N,M}(i)
N^{\sqrt{2} \gamma_{i-1} \beta}$.
Therefore on $B_{N,M,\delta}$ we get
\[
Z^{(\vec{\sigma}, \alpha)}_N(\beta) \ge\sum_{i=1}^{M}
\bigl(1-N^{\mathcal{E}(\gamma_{i})-\mathcal{E}(\gamma
_{i-1})+2\delta} \bigr) N^{P_{\beta}(\gamma_{i-1})-\delta}.
\]
This yields on $B_{N,M,\delta}$
\[
\label{eq:tm3-eq3} f^{(\vec{\sigma}, \alpha)}_N(\beta) \ge\frac
{\log(1-N^{ \max_{1 \le i
\le M}\{\mathcal{E}(\gamma_{i})-\mathcal{E}(\gamma_{i-1})\} + 2
\delta
})}{\log N}+
\max_{0 \le i \le M-1} P_{\beta}(\gamma_i)-\delta.
\]
Since for $\delta$ in (\ref{eq:tm3-choicedelta})
\[
\lim_{N \to\infty}(\log N)^{-1}\log\bigl(1-N^{ \max_{1 \le i \le
M}\{
\mathcal{E}(\gamma_{i})-\mathcal{E}(\gamma_{i-1})\} + 2 \delta}
\bigr)=0,
\]
the choices of $M$, $\delta$ in \eqref{eq:tm3-choiceM} and \eqref
{eq:tm3-choicedelta} give that
$f^{(\vec{\sigma}, \alpha)}_N(\beta) -f^{(\vec{\sigma}, \alpha
)}(\beta)
> -\nu$ on $B_{N,M,\delta}$ for $N$ large enough. Therefore, (\ref
{eq:tm3-eq1}) is a consequence of Lemma~\ref{lem:thm3-Bonenv}.

\textit{Proof of the upper bound} (\ref{eq:tm3-eq2}).
Observe first that the partition function $ Z^{(\vec{\sigma},
\alpha)}_N(\beta)$ satisfies on $B_{N,M,\delta}$
\[
Z^{(\vec{\sigma}, \alpha)}_N(\beta) \le\sum_{i=1}^{M}K_{N,M}(i)
N^{\sqrt{2} \gamma_{i} \beta}+ N,
\]
the second term coming from the negative values of the field.
Since $\mathcal E(0)=1$, on $B_{N,M,\delta}$ and for $N$ large enough,
we have using \eqref{eq:tm3-choicedelta}
\[
K_{N,M}(1)\geq N^{1-\delta}-N^{\mathcal E(\gamma_1)+\delta}\geq
\tfrac
{1}{2} N^{1-\delta},
\]
thus $N\leq2K_{N,M}(1)N^{\delta}$.
Moreover, on $B_{N,M,\delta}$ the random variable $K_{N,M}(i)$ are less
than $N^{\mathcal{E}(\gamma_{i-1})+\delta}$ for all $1 \le i \le M$.
The two last observations imply by the choice of $\delta$
\[
Z^{(\vec{\sigma}, \alpha)}_N(\beta) \le\sum_{i=1}^{M}
K_{N,M}(i) N^{\sqrt{2} \gamma_{i} \beta} + 2K_{N,M}(1)N^{\delta}
\le3 \sum_{i=1}^{M} N^{\mathcal{E}(\gamma_{i-1})+\sqrt{2} \gamma_{i}
\beta+\delta}.
\]
Therefore, on the event $B_{N,M,\delta}$, we get
\[
f^{(\vec{\sigma}, \alpha)}_N(\beta) \le\frac{\log(3M)}{\log N}+
\max
_{\gamma\in[0,\gamma_{\mathrm{max}} ]} P_{\beta}(\gamma) + \frac
{\sqrt{2} \beta}{M} +\delta.
\]
Recalling \eqref{eq:sup=free} and since $\lim_{N \to\infty}(\log
N)^{-1}\log(2M)=0$, the choices of $M$ and $\delta$ in (\ref
{eq:tm3-choiceM-bis}) and (\ref{eq:tm3-choicedelta}) imply that
$f^{(\vec{\sigma}, \alpha)}_N(\beta) - f^{(\vec{\sigma}, \alpha
)}(\beta
)< \nu$ on $B_{N,M,\delta}$ for $N$ large enough. Therefore (\ref
{eq:tm3-eq2}) is a consequence of Lemma \ref{lem:thm3-Bonenv}.

\begin{appendix}
\section*{Appendix}\label{app}

\subsection{Gaussian estimates, large deviation result and integration
by part}

%
\begin{lemmas}[(see, e.g., \cite{durrett})]
\label{lem: gaussian}
Let $X$ be a standard Gaussian random variable.
For any $a>0$, we have
\[
\frac{(1-2a^{-2})}{\sqrt{2 \pi} a} \ee^{-a^2/2} \le\p( X \ge a)
\le\frac{1}{\sqrt{2 \pi} a}
\ee^{-a^2/2}.
\]
\end{lemmas}

%
\begin{lemmas}[(see, e.g., \cite{bennett})]
\label{lem: LD}
Let $Z_1,\ldots, Z_n$ be i.i.d. real valued random variables satisfying
${\mathbb E}[ Z_i] =0$, $\sigma^2={\mathbb E}[ Z_i^2]$ and $\Vert Z_i
\Vert_\infty\le
1$. Then for any \mbox{$t>0$},
\[
\p\Biggl( \Biggl\llvert\sum_{i=1}^{n}
Z_i \Biggr\rrvert\ge t \Biggr) \le2 \exp\biggl\{ -
\frac{t^2}{2n\sigma^2 + 2t/3} \biggr\}.
\]
\end{lemmas}

%
\begin{lemmas}[(see, e.g., the Appendix of \cite{talagrand})]
\label{lem: GaussianIbP}
Let $(X,Z_1,\ldots, Z_d)$ be a centered Gaussian random vector. Then,
for any $C^1$ function $F\dvtx{\mathbb R}^d \mapsto{\mathbb R}$, of
moderate growth at
infinity, we have
\[
{\mathbb E}\bigl[ X F(Z_1,\ldots, Z_d) \bigr] = \sum
_{i=1}^d {\mathbb E}[ X Z_i ] {\mathbb E}\biggl[
\frac{\partial F}{\partial z_i}(Z_1,\ldots, Z_d) \biggr].
\]
\end{lemmas}

\subsection{Proof of Lemma \texorpdfstring{\protect\ref{lem:covariance-acc-perturbed}}{2.2}}
\label{proof-covariance-acc-perturbed}
Recall that $0<\varepsilon=1/N<1/2$, and $ \alpha\in(0,1)$.
Also by definition, $\Vert x' -x \Vert= \varepsilon^{q(x,x')}$.

It is clear that $\E[\tilde X_xX_x]=E[(\tilde X_x)^2]$, which is the
variance of the centered Gaussian random variable $\mu(A_{\varepsilon}
(x) \setminus A_{\varepsilon^{ \alpha}} (x))$. This variance can be
computed and equals
\[
\int_{\varepsilon}^{\varepsilon^{ \alpha}} y^{-1} \,\d y = [ \log y
]_{\varepsilon}^{\varepsilon^{ \alpha}}= (1-\alpha) \log N.
\]
For the covariance, observe that $\E[\tilde X_x X_{x'}]$ is equal to
the variance of the random variable $\mu((A_{\varepsilon} (x)
\setminus
A_{\varepsilon^{ \alpha}} (x)) \cap A_{\varepsilon} (x'))$.
If $ \varepsilon< \ell= \Vert x' -x \Vert< \varepsilon^{ \alpha}$
[i.e., $ \alpha< q(x,x') \le1$], then the subsets intersect in
between the lines $y=\varepsilon$ and $y=\varepsilon^{ \alpha}$, thus
\[
\E[\tilde X_x X_{x'}] =\int_{ \ell}^{\varepsilon^{ \alpha}}
\frac{y-\ell
}{y^2} \,\d y = [ \log y ]_{ \ell}^{\varepsilon^{ \alpha}}+ \ell
\biggl[ \frac{1}{y} \biggr]_{ \ell}^{\varepsilon^{ \alpha}}=
\bigl(q
\bigl(x,x'\bigr)- \alpha\bigr) \log N + O_N(1).
\]
Finally, if $ \ell= \Vert x' -x \Vert\ge\varepsilon^{ \alpha}$ [i.e.,
$0 \le q(x,x') \le\alpha$], then the set $(A_{\varepsilon} (x)
\setminus A_{\varepsilon^{ \alpha}} (x)) \cap A_{\varepsilon} (x')$ is
empty and thus $\E[\tilde X_x X_{x'}]=0$.

\subsection{A key property of the perturbed models}
The following lemma is a key tool to approximate the Gaussian field we
consider by a tree.
Indeed the difference between the contribution to the Gaussian field at
a certain scale for two points that are close
can be explicitly computed by integrating parallelograms (see Figure~\ref{fig:intersect-lemma} below) and is shown to be small.
%
%
\begin{lemmas}
\label{lem: lem12}
Fix $\alpha',\alpha''$ as in Lemma \ref{lem:recurrence}, $u$ such that
$\alpha' <u< \alpha''$ and $\delta\in(0,1)$.
Then for all $x,x' \in\mathcal{X}_{\varepsilon}$ such that $\Vert x-x'
\Vert\le\delta\varepsilon^u$, we have
\[
\operatorname{Var} \bigl( Y_{x}(u) - Y_{x'}(u) \bigr) \le2
\overline\sigma^2 \delta,
\]
where $\overline\sigma$ denotes an upper bound for the $\sigma_i$'s.
\end{lemmas}

\begin{pf} Writing $A:= A_{\varepsilon^u} (x) \Delta A_{\varepsilon
^u} (x')$, we have
\begin{eqnarray*}
\operatorname{Var} \bigl( Y_{x}(u) - Y_{x'}(u) \bigr) &\le&
\overline\sigma^2 \int_A y^{-2} \,\d s
\,\d y = 2 \overline\sigma^2\bigl \Vert x-x' \bigr\Vert\int
_{\varepsilon^u}^{\infty} y^{-2} \,\d y
\\
&=& 2 \overline\sigma^2 \frac{\Vert x-x' \Vert}{\varepsilon^u} \le
2 \overline
\sigma^2 \delta,
\end{eqnarray*}
which completes the proof of the lemma.
\end{pf}


%
\begin{figure}

\includegraphics{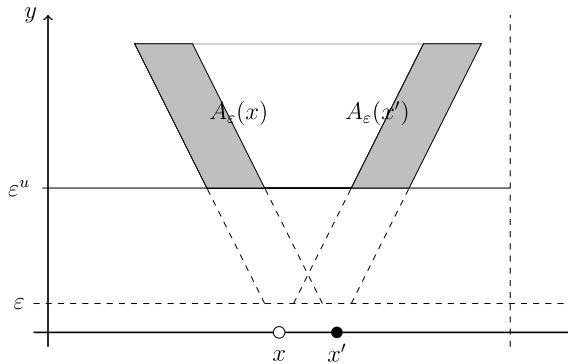}

\caption{The error terms in the tree approximation correspond to the
two grey parallelograms in Lemma \protect\ref{lem: lem12}.}
\label{fig:intersect-lemma}
\end{figure}
\end{appendix}

\section*{Acknowledgments} The authors thank Yan Fyodorov, Nicola
Kistler, Irina Kurkova and Vincent Vargas for helpful discussions.
O. Zindy would like to thank the Courant
Institute of Mathematical Science and the Universit\'e de Montr\'eal
for hospitality and financial support.

%



\printaddresses

\end{document}